\documentclass[a4paper,fleqn]{cas-sc}

\usepackage[authoryear,longnamesfirst]{natbib}

\usepackage{hyperref}
\hypersetup{
    colorlinks,%
    citecolor=black,%
    filecolor=black,%
    linkcolor=black,%
    urlcolor=black,
    breaklinks=true
} 
\usepackage{url}
\usepackage{enumerate}
\usepackage[shortlabels]{enumitem}
\usepackage{amsmath,amssymb}
\usepackage{amsfonts}
\usepackage[dvipsnames]{xcolor}
\graphicspath{{images/}}
\usepackage{bm}
\usepackage{soul}
\usepackage[percent]{overpic}
\usepackage{caption}
\usepackage{subcaption}
\usepackage{multirow}
\usepackage{scalerel}
\usepackage{mathtools}
\usepackage{import}
\usepackage{calc}
\usepackage{adjustbox}
\usepackage{listings}
\usepackage{algorithm}
\usepackage{algpseudocode}
\usepackage{tikz}
\usetikzlibrary{arrows,chains,scopes,positioning,matrix,shapes,calc}

\makeatletter
\tikzset{join/.code=\tikzset{after node path={%
\ifx\tikzchainprevious\pgfutil@empty\else(\tikzchainprevious)%
edge[every join]#1(\tikzchaincurrent)\fi}}}
\makeatother

\tikzset{>=stealth',every on chain/.append style={join},
         every join/.style={->}}

\usepackage{lscape}
\usepackage{rotating}
\usepackage{appendix}
\newcommand*{\colorboxed}{}
\def\colorboxed#1#{%
  \colorboxedAux{#1}%
}
\newcommand*{\colorboxedAux}[3]{%
  \begingroup
    \colorlet{cb@saved}{.}%
    \color#1{#2}%
    \boxed{%
      \color{cb@saved}%
      #3%
    }%
  \endgroup
}
\newcommand*\circled[1]{\tikz[baseline=(char.base)]{
            \node[shape=circle,draw,inner sep=1pt] (char) {#1};}}

\newcommand{\mbf}[1]{\mathbf{#1}} 

\newcommand{\what}[1]{\widehat{#1}}
\newcommand{\wtilde}[1]{\widetilde{#1}}

\newcommand{\Omegahat}{\what{\Omega}}
\newcommand{\bnabla}{\bm{\nabla}}

\newcommand{\ovbar}[1]{\overline{#1}}

\newcommand{\xhat}{\what{x}}
\newcommand{\yhat}{\what{y}}

\newcommand{\ipi}[3]{\int_{#3} #1 #2\,d\Omega}

\newcommand{\ipiv}[3]{\int_{#3} {#1} \cdot {#2}\,d\Omega}
\newcommand{\ipibdom}[3]{\int_{#3}#1#2\,d\Gamma}

\begin{document}

\shorttitle{An orthonormal and hierarchical divergence-free polynomial basis}    

\shortauthors{Sreevatsa Anantharamu and Krishnan Mahesh}

\title[mode = title]{Arnoldi-based orthonormal and hierarchical divergence-free
  polynomial basis and its applications}




\author[1]{Sreevatsa Anantharamu}


\address[1]{Aerospace Engineering and Mechanics, University of Minnesota - Twin Cities, 55455, USA}

\author[1]{Krishnan Mahesh}\cormark[1]


\cortext[cor1]{Corresponding author} 



\begin{abstract}
  This paper presents a methodology to construct a divergence-free polynomial
  basis of an arbitrary degree in a simplex (triangles in 2D and tetrahedra in
  3D) of arbitrary dimension. It allows for fast computation of all numerical
  solutions from degree zero to a specified degree $k$ for certain PDEs. The
  generated divergence-free basis is orthonormal, hierarchical, and robust in
  finite-precision arithmetic. At the core is an Arnoldi-based procedure. It
  constructs an orthonormal and hierarchical basis for multi-dimensional
  polynomials of degree less than or equal to $k$. The divergence-free basis is
  generated by combining these polynomial basis functions. An efficient
  implementation of the hybridized BDM mixed method is developed using these
  basis functions. Hierarchy allows for incremental construction of the global
  matrix and the global vector for all degrees (zero to $k$) using the local
  problem solution computed just for degree $k$. Orthonormality and
  divergence-free properties simplify the local problem. PDEs considered are
  Helmholtz, Laplace, and Poisson problems in smooth domains and in a corner
  domain. These advantages extend to other PDEs such as incompressible Stokes,
  incompressible Navier-Stokes, and Maxwell equations.
\end{abstract}

\maketitle

\section{Introduction}
Divergence-free vector fields occur in several problems. For e.g., the fluid
velocity field in an incompressible fluid flow, the solid velocity in an
incompressible solid deformation, the magnetic field around an electric current,
and the steady state heat flux in a conducting medium with no volumetric
sources. Approximating such vector fields using divergence-free basis functions
is advantageous. It reduces the number of global degrees of freedom while
computing approximate solutions to their partial differential equations (PDEs)
\citep{cockburn2004locally,cockburn2010comparison}. While interpolating
experimental measurements \citep{gesemann2016noisy,agarwal2021reconstructing},
it yields reconstructed fields that are consistent with the problem
physics. This paper discusses a procedure to construct a divergence-free
polynomial basis that is well-conditioned, orthonormal, and hierarchical for
arbitrary polynomial degree in a simplex.

Below is a simple exercise to construct a linear monomial divergence-free basis for two dimensions. Consider the two dimensional linear monomial basis:
\begin{equation*}
\Bigg\{\underbrace{\begin{pmatrix}1\\0\end{pmatrix}}_{(1)},\underbrace{\begin{pmatrix}0\\1\end{pmatrix}}_{(2)},\underbrace{\begin{pmatrix}x\\0\end{pmatrix}}_{(3)},\underbrace{\begin{pmatrix}y\\0\end{pmatrix}}_{(4)},\underbrace{\begin{pmatrix}0\\x\end{pmatrix}}_{(5)},\underbrace{\begin{pmatrix}0\\y\end{pmatrix}}_{(6)}\Bigg\}.
\end{equation*}
The basis functions (1), (2), (4), and (5) are divergence-free but (3) and (6) are not. However, basis functions (3) and (6) can be combined as 
\begin{equation*}
\begin{pmatrix}x\\0\end{pmatrix}-\begin{pmatrix}0\\y\end{pmatrix}=\begin{pmatrix}x\\-y\end{pmatrix}
\end{equation*}
to be divergence-free. This yields the following linear monomial divergence-free basis:
\begin{equation*}
\Bigg\{\underbrace{\begin{pmatrix}1\\0\end{pmatrix}}_{(1)},\underbrace{\begin{pmatrix}0\\1\end{pmatrix}}_{(2)},\underbrace{\begin{pmatrix}y\\0\end{pmatrix}}_{(3)},\underbrace{\begin{pmatrix}0\\x\end{pmatrix}}_{(4)},\underbrace{\begin{pmatrix}x\\-y\end{pmatrix}}_{(5)}\Bigg\}.
\end{equation*}
Similarly, the following quadratic monomial divergence-free basis:
\begin{equation*}
\Bigg\{\underbrace{\begin{pmatrix}1\\0\end{pmatrix}}_{(1)},\underbrace{\begin{pmatrix}0\\1\end{pmatrix}}_{(2)},\underbrace{\begin{pmatrix}y\\0\end{pmatrix}}_{(3)},\underbrace{\begin{pmatrix}0\\x\end{pmatrix}}_{(4)},\underbrace{\begin{pmatrix}x\\-y\end{pmatrix}}_{(5)},\underbrace{\begin{pmatrix}y^2\\0\end{pmatrix}}_{(6)},\underbrace{\begin{pmatrix}0\\x^2\end{pmatrix}}_{(7)},\underbrace{\begin{pmatrix}x^2/2\\-xy\end{pmatrix}}_{(8)},\underbrace{\begin{pmatrix}-xy\\y^2/2\end{pmatrix}}_{(9)}\Bigg\}
\end{equation*}
can be constructed from the two-dimensional quadratic monomial basis. This procedure can be generalized to arbitrary polynomial degree and spatial dimension; see the MATLAB program \texttt{mondivfreebf} in figure \ref{fig:mondivfreebf} of appendix \ref{app:mondivfreebf}. However, it turns out that the resulting monomial divergence-free basis functions perform poorly in finite precision arithmetic for high polynomial degrees.

\begin{figure}[t]
\centering
\begin{subfigure}[b]{0.3\linewidth}
\centering
\adjustbox{max width=\linewidth,trim=0cm 0cm 0cm 0cm,clip}{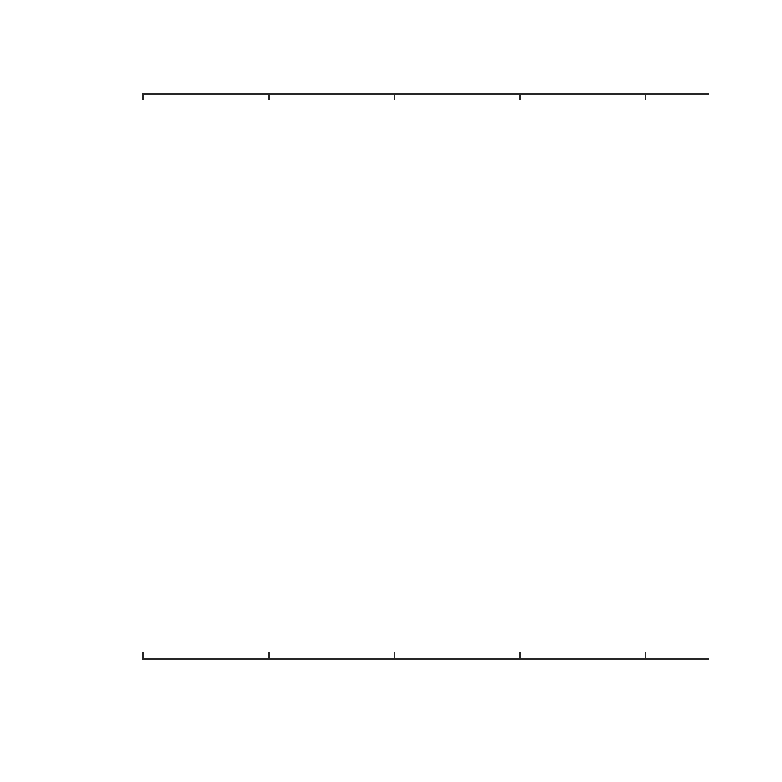}
\caption{}
\end{subfigure}\hfill
\begin{subfigure}[b]{0.3\linewidth}
\centering
\adjustbox{max width=\linewidth,trim=0cm 0cm 0cm 0cm,clip}{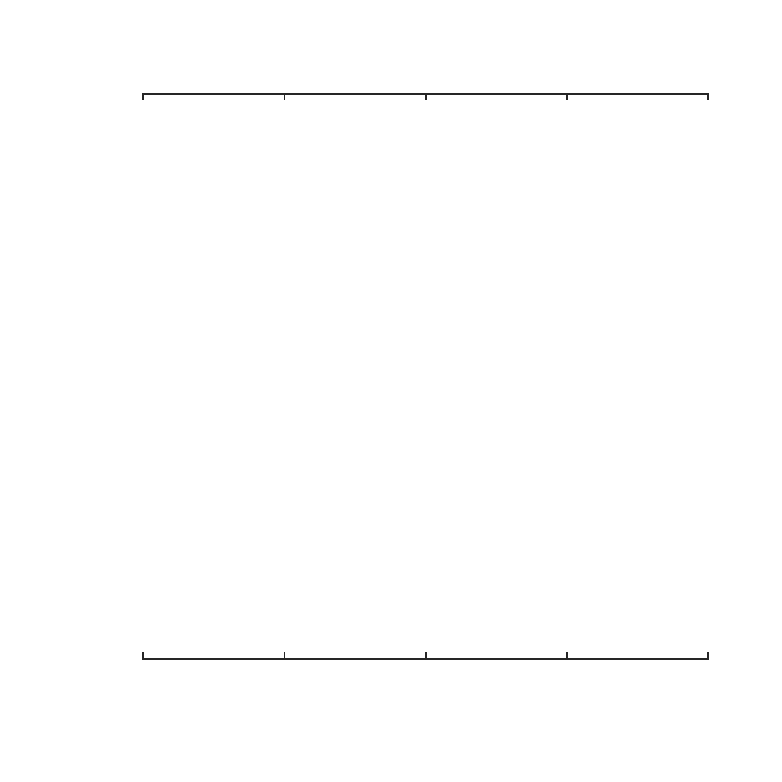}
\caption{}
\end{subfigure}\hfill
\begin{subfigure}[b]{0.35\linewidth}
\centering
\adjustbox{max width=\linewidth,trim=0cm 0.7cm 0.7cm 1cm,clip}{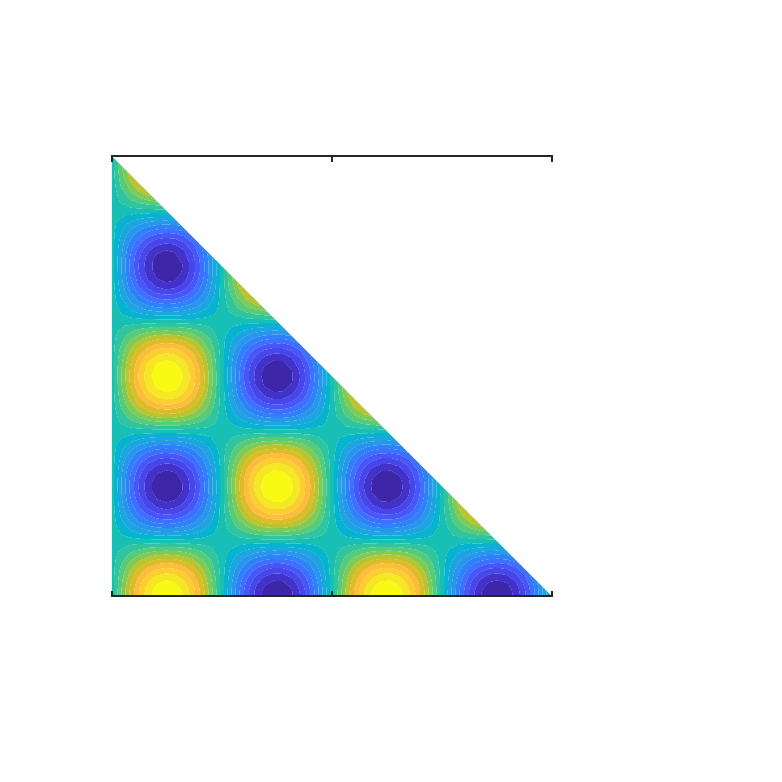}
\caption{}
\end{subfigure}
\caption{Projection of  $(\sin(4\pi x)\cos(4\pi y),-\cos(4\pi x)\sin(4\pi y))$ onto the space of divergence-free polynomials in the unit triangle.  (a) Error v/s polynomial degree for the monomial divergence-free basis. (a) Error v/s polynomial degree for the Arnoldi-based divergence-free basis.  (c) Contours of $x$-component of the projection computed with the Arnoldi-based divergence-free basis for polynomial degree 40.}
\label{fig:l2projdivfree}
\end{figure}

Consider the $L^2$ projection of the divergence-free function $\vec{g}=(\sin(4\pi x)\cos(4\pi y),-\cos(4\pi x)\sin(4\pi y))$ onto the space of divergence-free polynomials in the unit triangle defined by the nodes (0,0), (1,0), and (0,1). Figure \ref{fig:l2projdivfree}a shows the maximum error in the projection computed with the monomial divergence-free basis as a function of polynomial degree. The error decreases to around $10^{-3}$ for degree 20. After degree 20, only 1-3 significant digits of accuracy are obtained. This is because of finite precision error. The condition number of the mass matrix with the monomial divergence-free basis increases exponentially with degree, and therefore, the finite precision error also grows exponentially.

Suppose, instead of using the monomial divergence-free basis, we use the
proposed divergence-free basis. The projection error decreases all the way down
to machine precision; see figure \ref{fig:l2projdivfree}b. The contours of
$x$-component of the projection computed with this basis are shown in figure
\ref{fig:l2projdivfree}c for degree 40. This demonstrates the robustness of our
divergence-free basis in finite precision arithmetic. To build a
well-conditioned divergence-free basis, we combine orthonormal polynomial
\citep{gautschi2004orthogonal} basis functions. To generate these input
orthonormal polynomial basis functions in a simplex of arbitrary dimension, we
propose a simple Arnoldi-based procedure. This procedure is an extension of the
one-dimensional Arnoldi/Stieltjes process \citep{gautschi1982generating}
discussed in lecture 37 of \cite{trefethen1997numerical}.

\begin{figure}[t]
\centering
\begin{subfigure}[b]{0.3\linewidth}
\centering
\adjustbox{max width=\linewidth,trim=0cm 0cm 0cm 0cm,clip}{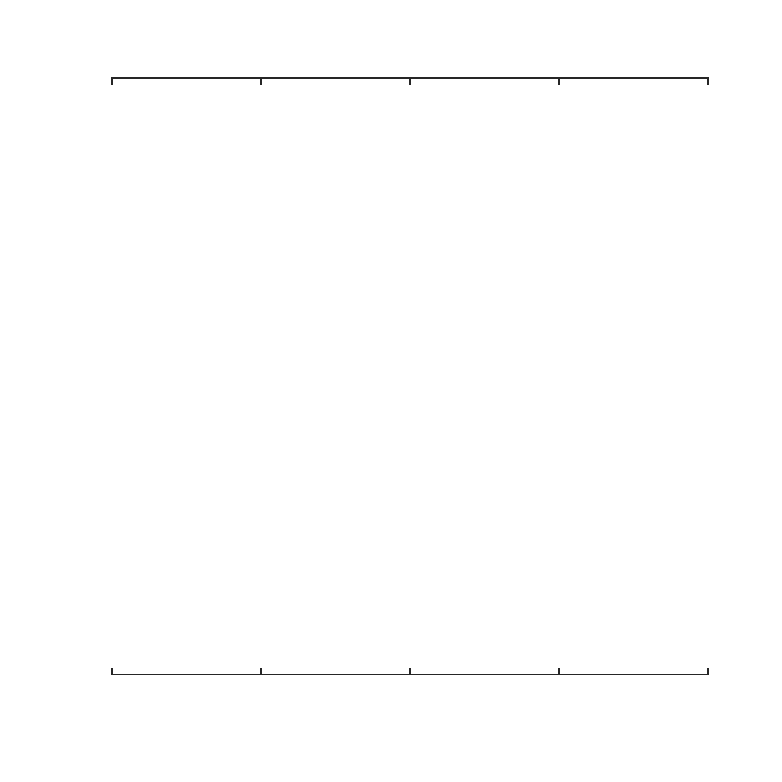}
\caption{}
\end{subfigure}\hfill
\begin{subfigure}[b]{0.3\linewidth}
\centering
\adjustbox{max width=\linewidth,trim=0cm 0cm 0cm 0cm,clip}{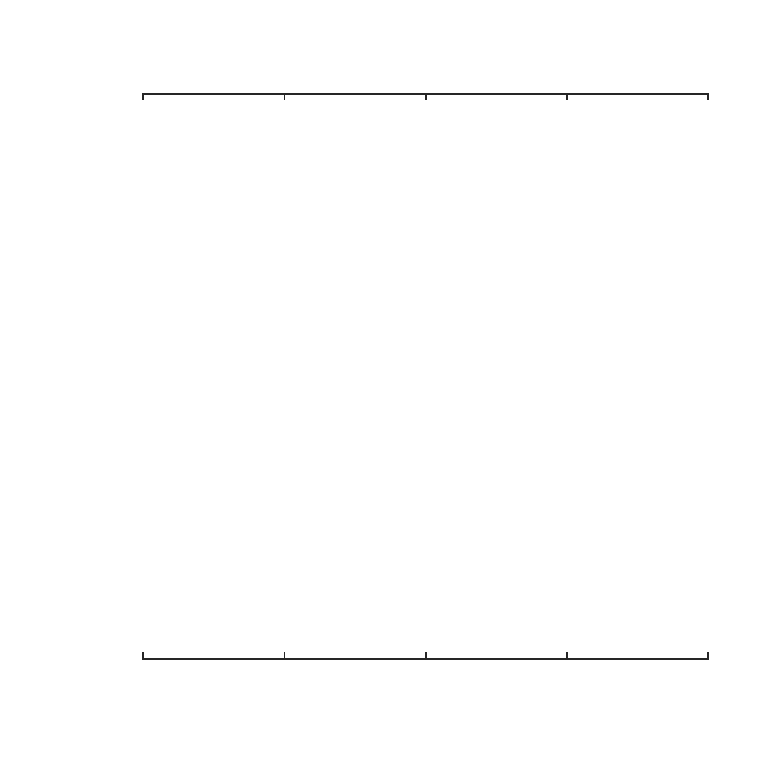}
\caption{}
\end{subfigure}
\begin{subfigure}[b]{0.35\linewidth}
\centering
\adjustbox{max width=\linewidth,trim=0cm 0.7cm 0.7cm 1cm,clip}{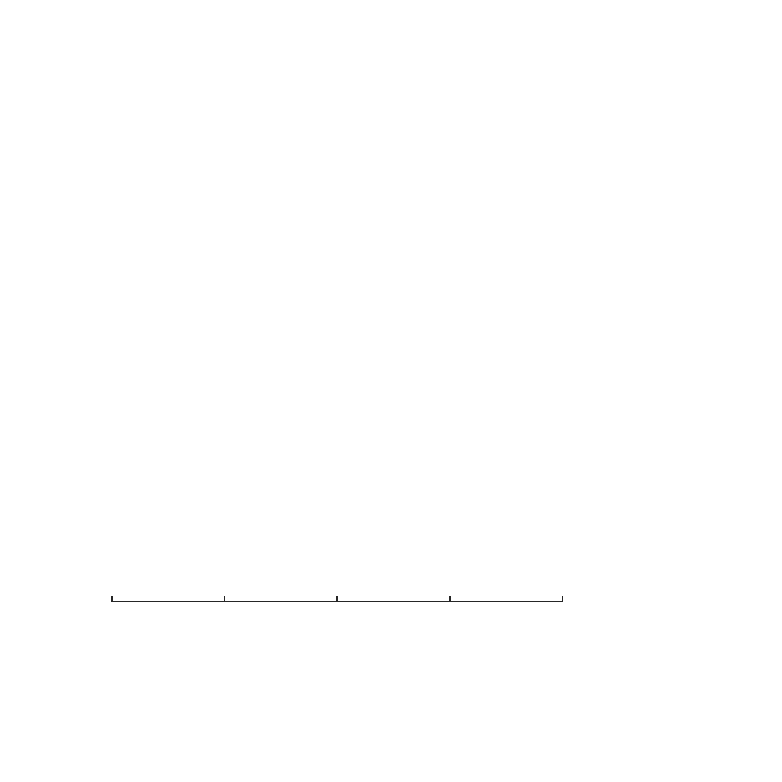}
\caption{}
\end{subfigure}\hfill
\caption{Laplace problem with corner singularity. (a) Mesh. (b) Error in the computed scalar at $(0.99,0.99)$ v/s polynomial degree. (c) Contours of $x$-component of the flux computed with polynomial degree eight. }
\label{fig:lap_corner_sing}
\end{figure}

The advantages of using our basis for numerical solutions of PDEs are
demonstrated for the Laplace problem with corner singularity. This problem is
taken from \cite{gopal2019new}. The domain is L-shaped (see figure
\ref{fig:lap_corner_sing}a). Dirichlet boundary conditions of $x^2$ are used on
all the boundaries, and the resulting solution has a singularity at the
re-entrant corner $(1,1)$. In $2019$, \cite{gopal2019new} called for finite
element method (FEM) solutions to this problem. Specifically, they asked for a
computation of the scalar at $(0.99,0.99)$ -- a point close to the re-entrant
corner. They report `\dots all respondents were able to calculate a solution to
two to four significant digits of accuracy, only two came close to eight
digits. For example, one researcher used 158,997 fifth-order triangular elements
near the re-entrant corner and achieved six correct digits \dots'.  Our results
computed by using the proposed divergence-free basis in the hybridized
Brezzi-Douglas-Marini (BDM) \citep{brezzi1985two} mixed finite element method
(FEM) are shown in figure \ref{fig:lap_corner_sing}. The approximation with
polynomial degree eight is accurate up to 12 significant digits at the point
$(0.99,0.99)$ and we use just 1000 elements (mesh shown in figure
\ref{fig:lap_corner_sing}a). It takes just around four seconds to compute all
approximations from polynomial degree zero to eight (all computations for this
paper are performed in MATLAB on a desktop workstation with Intel(R) Core(TM)
i7-8700 CPU @ 3.20GHz and six cores).

We can compute all approximations in such a short time because our basis is
hierarchical and orthonormal. Hierarchy allows us to solve the local problem in
the hybridized BDM method just for polynomial degree eight and use its solution
to incrementally construct the element (and global) matrices and vectors for all
polynomial degrees from zero to eight. Orthonormality simplifies the local
problem solutions to just inner products instead of requiring matrix
inversions. Therefore, our computation is fast. Furthermore, the results
demonstrate exponential convergence near the singularity. \cite{gopal2019new}
note that hp-adaptive FEM can achieve exponential convergence near singularities
but requires advanced implementations. We, on the other hand, do not require
such advanced implementation. These advantages of our basis extend to several
other PDEs and to other hybridized FEM methods.


A word on the one-dimensional Arnoldi/Stieltjes process discussed in
\cite{trefethen1997numerical}. The development of this process begins by
recognizing that the space of one-dimensional polynomials of degree less than or
equal to $k$ given by span$\{1,x,x^2,...,x^k\}$ is a Krylov subspace. Define the
coordinate operator $\xhat$ as the operator that maps a one-dimensional
polynomial $f$ to another polynomial $xf$. The one-dimensional polynomial space
can then be rewritten as the Krylov subspace
span$\{1,\xhat 1,\xhat^2 1,...,\xhat^k 1\}$ generated by the linear operator
$\xhat$ and the starting polynomial `1'. The one-dimensional Arnoldi/Steiljets
process to generate an orthonormal basis for this space is:
\begin{algorithmic}[1]
\State $q_1=1$
\For{$j=1,\dots,k$}
\State $v=\xhat q_j$
\For{$i=1,\dots,j$}
\State $h_{i,j}=\int_0^1q_i^*v\,dx$
\State $v=v-h_{i,j}q_i$
\EndFor
\State $h_{j+1,j}=\left(\int_0^1|v|^2\,dx\right)^{1/2}$
\State $q_{j+1}=v/h_{j+1,j}$
\EndFor
\end{algorithmic}
\noindent The generated polynomials $q_1,\dots,q_{k+1}$ are a basis orthonormal
in the $L^2$ inner-product for polynomials of degree less than or equal to $k$
in the interval $[0,1]$. Since the operator $\xhat$ is hermitian, the above
Arnoldi process can be simplified to a Lanczos process. However, such
simplifications are not performed usually for numerical stability reasons. We
note that this process is rarely used to generate the one-dimensional
polynomials orthonormal in the $L^2$ inner-product, i.e., the Legendre
polynomials. Instead, analytical expressions for the coefficients $h_{i,j}$,
also called the recurrence relations, are used. Nevertheless, it still is a
powerful technique to generate one-dimensional polynomials orthonormal for an
arbitrary weighted inner-product.

The situation is a little different in a simplex of dimension larger than
one. Recurrence relations do exist for the coefficients analogous to $h_{i,j}$
to construct orthonormal polynomials in triangles and tetrahedra
\citep{olver2020fast,sherwin1995new,dubiner1991spectral}. However, these are
complicated compared to the one dimensional relations. Our Arnoldi-based process
is a simple alternative implementable in just a few lines of MATLAB code. An
advantage of this process is that it can construct an orthonormal basis not just
for the $L^2$ inner-product but for arbitrary weighted inner-product and for
even discrete inner-products. Another advantage is that it extends to a simplex
of arbitrary dimension.

The proposed Arnoldi-based process can also be seen as an extension of the
`Vandermonde with Arnoldi' idea of
\cite{brubeck2021vandermonde}. \cite{brubeck2021vandermonde} considered the
one-dimensional Vandermonde matrix problem and showed that despite using the
well-conditioned Chebyshev points to construct the Vandermonde matrix, the
computed approximation at high polynomial degrees can be significantly
contaminated by round-off error. To remedy this issue, they proposed an
Arnoldi-based procedure to solve the Vandermonde matrix problem. Using this
procedure, they were able to compute approximations that were accurate up to
machine precision. The $L^2$ projection numerical experiment whose results are
displayed in figure \ref{fig:l2projdivfree} demonstrate a similar accuracy in
finite precision arithmetic for the proposed Arnoldi-based process.

Finally, we note that \cite{ainsworth2018bernstein} constructed a divergence-free basis for triangles and tetrahedra using Bernstein polynomials. However, their construction yields a non-orthogonal basis, thus requiring matrix inversions in local problem solution in hybridized FEM methods while ours requires only inner products. Furthermore, their divergence-free basis is not hierarchical. Therefore, the local problem solution, construction of element matrices and vectors need to be computed separately for each polynomial degree, while we exploit the basis hierarchy to compute them efficiently.

The rest of the paper is organized as follows. In section
\ref{sec:divfreebasesconst}, we present the proposed method. Its numerical
implementation is given in section \ref{sec:impl}. Some remarks on the proposed
method are made in section \ref{sec:remarks}. Section \ref{sec:appln}
demonstrates some applications of the proposed basis. The paper is summarized in
section \ref{sec:summary}.


\section{Divergence-free polynomial basis construction} \label{sec:divfreebasesconst}



We are given $N_{ele}$ simplices (also referred to as `elements') in $d$
dimensions. The node coordinate matrix (a matrix whose rows store the coordinate
vector of the nodes of the element) of each element $e$ is $X^{(e)}$. We need to
construct an orthonormal and hierarchical basis for divergence-free polynomials
of degree less than or equal to some prescribed degree $k$ in each of these
elements. The required basis is the set of vector-valued polynomials
$\bm{\varphi}_1^{(e)},\dots,\bm{\varphi}_n^{(e)}$ in each element $e$. By an
orthonormal basis, we mean that any two divergence-free basis functions
$\bm{\varphi}_i^{(e)}$ and $\bm{\varphi}_j^{(e)}$ of element $e$ satisfy the
orthonormality relation
$\int_{\Omega^{(e)}}\bm{\varphi}^{(e)}_{i}\cdot\bm{\varphi}^{(e)}_{j}\,d\Omega=\delta_{ij}|\Omega^{(e)}|$,
where $\Omega^{(e)}$ denotes the $e^{th}$ element and $|\Omega^{(e)}|$ is its
volume. By a hierarchical basis, we mean that the basis functions are generated
incrementally for each degree up to $k$. A formal definition of this is that the
first $dC_d^{j+d}-C_d^{j-1+d}$ basis functions form a basis for divergence-free
polynomials of degree less than or equal to $j$, where $j$ is any degree
satisfying $0\leq j\leq k$. $n$ is the dimension of the basis in each
element. It equals $dC_d^{k+d}-C_d^{k-1+d}$ because the dimension of the set of
vector-valued polynomials is $dC_d^{k+d}$ and the divergence-free requirement
imposes $C_d^{k-1+d}$ constraints.

The proposed method can be summarized as follows:
\begin{enumerate}[label={\bfseries Step \arabic*:}]
\item Construct the divergence-free basis first in the reference element
  $\Omegahat$. This basis is the set of vector-valued polynomials
  $\bm{\varphi}_1,\dots,\bm{\varphi}_n$, where $\bm{\varphi}_i$ is the $i^{th}$
  basis function. The reference element is a unit simplex in $d$ dimensions. Its
  node-coordinate matrix is [zeros$(1,d)$; eye$(d)$] (in MATLAB notation). The
  goal is to combine polynomial basis functions that are orthonormal in the
  reference element to generate the divergence-free polynomial basis functions
  that are also orthonormal in the reference element. The orthonormal polynomial
  basis is the set of polynomials $q_1,\dots,q_p$, where $p=C_d^{k+d}$ and $q_i$
  is the $i^{th}$ basis function. To ensure that the constructed divergence-free
  basis functions are hierarchical, this idea is recursively applied for each
  polynomial degree up to $k$. The required polynomial basis functions ($q_i$)
  are constructed for each degree using an Arnoldi-based procedure. For degree
  $j=0$, there is only one polynomial basis function and it is $q_1=1$. There
  are three degree zero divergence-free polynomial basis functions given by
  $\bm{\varphi}_i=q_1\mbf{e}_i$, where $i=1,\dots,d$. For each degree
  $j=1,\dots,k$, do \textbf{Step} \textbf{1.1} and \textbf{1.2}.
  \begin{enumerate}[label={\bfseries Step 1.\arabic*:}]
  \item Compute the polynomial basis functions $q_{p_{j-1}+1},\dots,q_{p_j}$
    of degree $j$ using the previously computed basis functions
    $q_1,\dots,q_{p_{j-1}}$ and the Arnoldi-based procedure given below:
      \begin{algorithmic}[1]
        \State Set $c=p_{j-1}$
        \For{$i=1,\dots,d$}
          \State $j''=C_{d-i}^{j-1+d-i}$
          \For{$j'=1,\dots,j''$}
            \State $v=\xhat_iq_{p_{j-1}-j''+j'}$ \Comment{Equivalent of $Aq_j$ in the standard Arnoldi}
            \For{$j'''=1,\dots,c$} \Comment{Orthogonalization}
              \State $H_{j''',c}=(\int_{\Omegahat}q_{j'''}^*v\,d\Omegahat)d!$
              \State $v=v-H_{j''',c}q_{j'''}$
            \EndFor
            \State $H_{c+1,c}=(\int_{\Omegahat}|v|^2\,\,d\Omegahat)^{1/2}{d!}^{1/2}$
            \State $q_{c+1}=v/H_{c+1,c}$\Comment{Normalization}
            \State $c=c+1$
          \EndFor
        \EndFor
      \end{algorithmic}
      Here, $p_j=C_d^{j+d}$. $\xhat_i$ is the coordinate operator along the
      $i^{th}$ direction. It maps a polynomial $f$ to another polynomial
      $x_if$. The new set of polynomials $q_1,\dots,q_{p_j}$ form a basis for
      polynomials of degree less than or equal to $j$ (see appendix
      \ref{app:proofarpoly} for more discussion). The Gram-Schmidt procedure
      enforces them to be orthonormal to each other. They satisfy the
      orthonormality relation
      $\int_{\Omegahat}q_iq_jd\Omegahat=\delta_{ij}1/d!$, where $1/d!$ is the
      volume of the reference element. The computed $H_{i,j}$s are stored in a
      matrix $H$ of size $p\times (p-1)$. Similar to the standard Arnoldi
      \citep{saad2011numerical}, $H$ is an upper Hessenberg matrix, i.e.,
      $H_{i,j}=0$ for $i>j+1$.
    \item Compute the divergence-free basis functions
      $\bm{\varphi}_{n_{j-1}+1},\dots,\bm{\varphi}_{n_j}$ of degree $j$ (where
      $n_j=dC_d^{j+d}-C_d^{j-1+d}$) by combining the polynomial basis functions
      $q_1,\dots,q_{p_{j}}$ and using the previously computed divergence-free
      basis functions $\bm{\varphi}_1,\dots,\bm{\varphi}_{n_{j-1}}$ as follows:
    \begin{algorithmic}[1]
      \State Expand each $\bm{\varphi}_{\ell}$ as $\bm{\varphi}_{\ell}=\sum_{i=1}^d\sum_{r=1}^{p_j}N_{(i-1)p+r,\ell}q_r\bm{e}_i,$ for $\ell=n_{j-1}+1,\dots,n_j$. Here, $N_{(i-1)p+r,\ell}$ are the coefficients stored in a matrix $N$ of size $dp\times n$.
      \State Compute the coefficients $N_{(i-1)p+r,\ell}$ such that the set of functions $\bm{\varphi}_{n_{j-1}+1},\dots,\bm{\varphi}_{n_{j}}$ are:
      \begin{enumerate}
      \item divergence-free, i.e., $\int_{\Omegahat}q_i\bnabla\cdot\bm{\varphi}_{\ell}d\Omegahat=0$, for $\ell=n_{j-1}+1,\dots,n_j$ and $i=1,\dots,p_{j-1}$,
      \item orthonormal, i.e., $\int_{\Omegahat}\bm{\varphi}_{\ell}\cdot \bm{\varphi}_{\ell'}\,d\Omegahat=\delta_{\ell,\ell'}1/d!$, where $\ell=n_{j-1}+1,\dots,n_j$, $\ell'=1,\dots,n_j$, and 
      \item linearly independent.
      \end{enumerate}
    \end{algorithmic}
    The new set of vector-valued polynomials
    $\bm{\varphi}_1,\dots,\bm{\varphi}_{n_j}$ form a orthonormal basis for the
    divergence-free polynomials of degree less than or equal to $j$. They
    satisfy the orthonormality relation
    $\int_{\Omegahat}\bm{\varphi}_i\cdot\bm{\varphi}_jd\Omegahat=\delta_{ij}1/d!$. See
    appendix \ref{app:divfreestep1p2} for a discussion on this.
\end{enumerate}
  
\item Construct the divergence-free basis in each element $e$ using the
  divergence-free basis in the reference element and the node coordinate matrix
  $X^{(e)}$ of the element. The basis in element $e$ is given by the set of
  vector-valued polynomials
  $\bm{\varphi}^{(e)}_1,\dots,\bm{\varphi}^{(e)}_n$. To ensure that the
  constructed divergence-free basis is hierarchical, the basis functions are
  constructed incrementally for each degree up to $k$. For each element
  $e=1,\dots,N_{ele}$ and for each degree $j=0,\dots,k$, do \textbf{Step 2.1}.
  \begin{enumerate}[label={\bfseries Step 2.\arabic*:}]
  \item Construct the divergence-free basis functions
    $\bm{\varphi}^{(e)}_{n_{j-1}+1},\dots,\bm{\varphi}^{(e)}_{n_j}$ of degree
    $j$ using the previously computed basis functions
    $\bm{\varphi}^{(e)}_1,\dots,\bm{\varphi}^{(e)}_{n_{j-1}}$ and the reference
    element basis functions
    $\bm{\varphi}_{n_{j-1}+1},\dots,\bm{\varphi}_{n_j}$ as follows:
    \begin{algorithmic}[1]
      \State Expand each $\bm{\varphi}^{(e)}_{\ell}$ as $\bm{\varphi}^{(e)}_{\ell}=\sum_{i=1}^d\sum_{r=1}^{p_j}N^{(e)}_{(i-1)p+r,\ell}q_r(\bm{x}(\bm{x}^{(e)}))\bm{e}_i,$ for $\ell=n_{j-1}+1,\dots,n_j$. Here, $N^{(e)}_{(i-1)p+r,\ell}$ are the coefficients stored in a matrix $N^{(e)}$ of size $dp\times n$, and $\bm{x}(\bm{x}^{(e)})$ is the mapping that maps the element coordinates $\bm{x}^{(e)}$ to the reference coordinates.
      \State Compute the coefficients $N^{(e)}_{(i-1)p+r,\ell}$ using the reference element coefficients $N_{(i-1)p+r,\ell}$ and the node-coordinate matrix $X^{(e)}$ such that the set of functions $\bm{\varphi}^{(e)}_{n_{j-1}+1},\dots,\bm{\varphi}^{(e)}_{n_{j}}$ are:
      \begin{enumerate}
      \item divergence-free, i.e., $\int_{\Omega^{(e)}}q_i\bnabla^{(e)}\cdot\bm{\varphi}_{\ell}^{(e)}d\Omega=0$ for $\ell=n_{j-1}+1,\dots,n_j$ where, $\bnabla^{(e)}$ is the gradient in the element coordinates,
      \item orthonormal, i.e., $\int_{\Omega^{(e)}}\bm{\varphi}^{(e)}_{\ell}\cdot \bm{\varphi}^{(e)}_{\ell'}\,d\Omega=\delta_{\ell,\ell'}|\Omega^{(e)}|$, where $\ell=n_{j-1}+1,\dots,n_j$, $\ell'=1,\dots,n_j$, and 
      \item linearly independent.
      \end{enumerate}
    \end{algorithmic}
    The set of vector-valued polynomials
    $\bm{\varphi}^{(e)}_1,\dots,\bm{\varphi}^{(e)}_{n_j}$ form an orthonormal
    basis for divergence-free polynomials of degree less than or equal to
    $j$. They satisfy the orthonormality relation
    $\int_{\Omega^{(e)}}\bm{\varphi}^{(e)}_i\cdot\bm{\varphi}^{(e)}_jd\Omega=\delta_{ij}|\Omega^{(e)}|$. Note
    that the above step is very similar to \textbf{step 1.2}. Instead of using
    reference element quantities, we now use the quantities of element
    $e$. Therefore, the discussion in appendix \ref{app:divfreestep1p2} applies
    to this step but for element $e$ instead of the reference element. However,
    the implementation of these two steps (discussed in the next section) differ
    significantly.

  \end{enumerate}
\end{enumerate}

The implementation details of each step follows. We use the MATLAB notation to
simplify the discussion. MATLAB implementations of step 1 and 2 are also given
and discussed in appendices \ref{app:matlabstep1} and \ref{app:matlabstep2},
respectively.

\section{Implementation} \label{sec:impl}

The basis construction procedure in the above form requires symbolic
manipulation of polynomials. To instead use just arithmetic computation, some
modifications are made. Instead of symbolically storing the polynomials, its
value at the quadrature points inside the unit simplex is stored as a vector.
$q_i$ is now a vector that stores the value of the orthonormal polynomial at the
quadrature points. Similarly, $\bm{\varphi}_i$ is now a vector that stores the
value of each component of the divergence-free basis function at the quadrature
points one below the other.

Denote the coordinates and weights of the quadrature rule by $x(:,:)$ and
$w(:)$, respectively. $x(:.i)$ stores the $i^{th}$ component of the coordinate
vector of the quadrature points. These points and weights are generated inside
the simplex using the Duffy transformation \citep{duffy1982quadrature}
(described in appendix \ref{app:quadrule}). To ensure exact integration of
polynomials of degree $2k$ that occur in the integrand of the norms and
inner-products, $k+1$ points are chosen in each direction. The total number of
quadrature points is $(k+1)^d$. Thus, we have the following relations. The
$j^{th}$ component of the vector $q_i$ stores the value of the $i^{th}$
orthonormal polynomial at the $j^{th}$ quadrature point. The
$((\ell-1)(k+1)^d+j)^{th}$ component of the vector $\bm{\varphi}_i$ stores the
value of the $\ell^{th}$ component of the $i^{th}$ divergence-free basis
function at the $j^{th}$ quadrature point.

\subsection{\textbf{Step 1}}

For degree zero, the only polynomial basis vector $q_1$ is $ones((k+1)^d,1)$,
where ${ones}(m,n)$ denotes a matrix of size $m\times n$ storing the number
one. The degree zero divergence-free basis vectors in the reference element are
$\bm{\varphi}_i={kron}(\mbf{e}_i,q_1)$ for $i=1,\dots,d$, where ${kron}$ denotes
the Kronecker tensor product of the two input matrices. Allocate space for the
upper Hessenberg matrix $H$ (size $p\times(p-1)$), divergence-free constraint
matrix $C$ (size $p_{k-1}\times dp$), and the coefficient matrix $N$ (size
($dp\times n$)). Initialize the first $d$ columns of $N$ as
$N_{1:p:dp,1:d}=eye(d)$ to be consistent with the initialization of the first
$d$ divergence-free basis functions $\bm{\varphi}_1,\dots,\bm{\varphi}_d$. Here,
$eye(d)$ denotes the $d\times d$ identity matrix. Initialize the upper
Hessenberg matrix and the constraint matrix $C$ to zero. For each degree
$j=1,\dots,k$, do \textbf{step 1.1} and \textbf{1.2}.

\subsubsection{\textbf{Step 1.1}}
Generate the new polynomial basis vectors $q_{p_{j-1}},\dots,q_{p_j}$ using the
below Arnoldi-based procedure.
\begin{algorithmic}[1]
  \State Set $c=p_{j-1}$
  \For{$i=1,\dots,d$}
    \State $j''=C_{d-i}^{j-1+d-i}$
    \For{$j'=1,\dots,j''$}
      \State $v={diag}(x(:,i))q_{p_{j-1}-j''+j'}$\Comment{Compute the matrix-vector product}
      \For{$j'''=1,\dots,c$}\Comment{First round of orthogonalization}
        \State $H_{j''',c}=(q_{j'''}^H{diag}(w(:))v))^{1/2}d!$
        \State $v=v-H_{j''',c}q_{j'''}$
      \EndFor
      \For{$j'''=1,\dots,c$}\Comment{Second round of orthogonalization}
        \State $t=(q_{j'''}^H{diag}(w(:))v)^{1/2}d!$
        \State $H_{j''',c}=H_{j''',c}+t$
        \State $v=v-tq_{j'''}$
      \EndFor
      \State $H_{c+1,c}=(v^H{diag}(w(:))v))^{1/2}{d!}^{1/2}$
      \State $q_{c+1}=v/H_{c+1,c}$
      \State $c=c+1$
    \EndFor
  \EndFor
\end{algorithmic}
      
The vectors $q_1,\dots,q_{p_j}$ are the value of the first $p_j$ orthonormal
polynomials at the quadrature points. They satisfy the discrete orthonormality
relation $q_i^H{diag}(w(:))q_j=\delta_{ij}1/d!$. The coordinate operators
$\xhat_1,\dots,\xhat_d$ are replaced by their discrete equivalent which are the
diagonal matrices ${diag}(x(:,1))$, \dots, ${diag}(x(:,d))$. The continuous
$L^2$ inner products are replaced by their discrete equivalent which is the
weighted $\ell^2$ inner-product with ${diag}(w(:)$ as the weight matrix. To
orthogonalize the $q_j$s, we have used the modified Gram-Schmidt kernel with
reorthogonalization. The symbolic Arnoldi procedure discussed in the previous
section used the modified Gram-Schmidt kernel with no reorthogonalization. Both
are equivalent in exact arithmetic. In finite-precision arithmetic, the former
leads to vectors that are orthogonal up to machine precision while the latter
leads to vectors that are orthogonal up to machine epsilon $\times$ the
condition number of the upper Hessenberg matrix $H$
\citep{saad2011numerical}. Therefore, we use the former for better numerical
accuracy.

\subsubsection{\textbf{Step 1.2}} \label{subsubsec:1p2}

Generate the new divergence-free basis vectors
$\bm{\varphi}_{n_{j-1}+1},\dots,\bm{\varphi}_{n_j}$ as follows:
\begin{algorithmic}[1]
  \State Compute the index vector $ii=[]$ and
  \For{$i=1,\dots,d$}
    \State $ii=[ii,\,(i-1)p+1:(i-1)p+p_j]$
  \EndFor
  \For{$i=1,\dots,d$} \Comment{Compute the new columns of the constraint matrix}
    \For{$r=p_{j-1}+1,\dots,p_j$}
      \State $v=DX_iq_r$
      \For{$j'=1,\dots,p_{j-1}$}\Comment{First round of projection}
        \State $C_{j',(i-1)p+r}=q_{j'}^H{diag}(w(:))v$
        \State $v=v-C_{j',(i-1)p+r}q_{j'}$
      \EndFor
      \For{$j'=1,\dots,p_{j-1}$}\Comment{Second round of projection}
        \State $t=q_{j'}^H{diag}(w(:))v$
        \State $C_{j',(i-1)p+r}=C_{j',(i-1)p+r}+t$
        \State $v=v-tq_{j'}$
      \EndFor      
    \EndFor
  \EndFor
  \State Compute the coefficients by setting $N_{ii,n_{j-1}+1:n_j}$ to $null\left(\begin{bmatrix}C_{1:p_{j-1},ii}\\ N^T_{ii,1:n_{j-1}}\end{bmatrix}\right)$
  \For{$\ell=n_{j-1}+1,\dots,n_j$}\Comment{Compute the new basis vectors}
    \State $\bm{\varphi}_{\ell}$ = $\sum_{i=1}^d\sum_{r=1}^{p_j}N_{(i-1)p+r,\ell}{kron}(\mbf{e}_i,q_r)$
  \EndFor
\end{algorithmic}
The vectors $\bm{\varphi}_1,\dots,\bm{\varphi}_{n_j}$ are the value of the first
$n_j$ divergence-free basis functions at the quadrature points. They satisfy the
discrete orthonormality relation
$\bm{\varphi}_i^H{diag}(w_d(:))\bm{\varphi}_j=\delta_{ij}1/d!$, where
$w_d={kron}({ones}(d,1),w)$. See appendix \ref{app:1p2} for a derivation of the
above algorithm from that in \textbf{step 1.2} of section
\ref{sec:divfreebasesconst}. Here, $DX_i$ is the derivative matrix along the
$i^{th}$ coordinate direction. Multiplying it with the polynomial basis vector
$q_j$ yields the value of the partial derivative of the $j^{th}$ orthonormal
polynomial along the $i^{th}$ direction at the quadrature points. Its
construction is described in appendix \ref{app:dermat}. The function $null()$
returns an orthonormal basis for the null-space of the input matrix. Note that
to compute the entries of the divergence-free constraint matrix $C$, we have
used the modified Gram-Schmidt kernel with reorthogonalization for better
numerical accuracy.

\subsection{\textbf{Step 2}}
For each element $e=1,\dots,N_{ele}$, allocate space for the coefficient matrix
$N^{(e)}$ (size $dp\times n$) and initialize it to zero. Compute the Jacobian
matrix $F^{(e)}$ (size $d\times d$). Its entries are
$F^{(e)}_{i,j}=\partial x_i/\partial x_j^{(e)}$. For each degree $j=0,\dots,k$,
do \textbf{step 2.1}.

\subsubsection{\textbf{Step 2.1}}\label{subsubsec:2p1}
Generate the new divergence-free basis vectors
$\bm{\varphi}^{(e)}_{n_{j-1}+1},\dots,\bm{\varphi}^{(e)}_{n_j}$ of element $e$
as follows.
\begin{algorithmic}[1]
  \State Initialize $\ovbar{N}$ to a zero matrix of size $dp_j\times (n_j-n_{j-1})$
  \For{$j=1,\dots,d$}\Comment{{Compute linear combination of the rows of $N$}}
    \For{$i=1,\dots,d$}
      \State $\ovbar{N}_{(i-1)p_j+1:ip_j,:}=\ovbar{N}_{(i-1)p_j+1:ip_j,:}+F^{(e)}_{i,j}N_{(j-1)p+1:(j-1)p+p_j,n_{j-1}+1:n_j}$
    \EndFor
  \EndFor
  \State Compute the index vector as $ii=[]$ and
  \For{$i=1,\dots,d$}
    \State $ii=[ii\hspace{.5em}(i-1)p+1:(i-1)p+p_j]$
  \EndFor
  \State Project $\ovbar{N}$ to be orthogonal to the previous columns of $N^{(e)}$ by computing $\ovbar{N}=\ovbar{N}-N^{(e)}_{ii,1:n_{j-1}}\left(N^{{(e)}^H}_{ii,1:n_{j-1}}\ovbar{N}\right)$
  \State Orthonormalize the columns of $\ovbar{N}$ by computing $\ovbar{N}=orth(\ovbar{N})$.
  \State Check for error in orthogonality $\ovbar{N}$ by computing $T=\ovbar{N}^HN^{{(e)}^H}_{ii,1:n_{j-1}}$ and $tt=norm(T,\text{`Inf'})$.
  \If{$tt>10^{-13}$}
    \State Reorthogonalize by computing $\ovbar{N}=\ovbar{N}-N^{(e)}_{ii,1:n_{j-1}}\left(N^{{(e)}^H}_{ii,1:n_{j-1}}\ovbar{N}\right)$ and setting $\ovbar{N}=orth(\ovbar{N})$.
  \EndIf
  \For{$\ell=n_{j-1}+1,\dots,n_j$}\Comment{Compute the new basis vectors}
    \State $\bm{\varphi}^{(e)}_{\ell}$ = $\sum_{i=1}^d\sum_{r=1}^{p_j}N^{(e)}_{(i-1)p+r,\ell}{kron}(\mbf{e}_i,q_r)$
  \EndFor
\end{algorithmic}
Here, $orth()$ returns an orthonormal basis for the span of the input
matrix. See appendix \ref{app:2p1} for a derivation of the above algorithm from
the algorithm in \textbf{step 2.1} of section \ref{sec:divfreebasesconst}. To
orthonormalize the new columns of the coefficient matrix against its previous
columns, a classical Gram-Schmidt-type kernel is used. We specifically chose
classical Gram-Schmidt instead of modified Gram-Schmidt because the former is
faster than the latter even though both have the same operation
acount. Classical Gram-Schmidt is faster because it predominantly uses
matrix-matrix multiplications which yield higher FLOPS compared to the modified
Gram-Schmidt which mainly uses matrix-vector multiplications. This difference in
performance was found to be crucial because the above algorithm is executed for
each element in the mesh. Using a modified Gram-Schmidt-type kernel led to a
slow down of factor 10 in some cases. But an issue with classical Gram-Schmidt
is that it can lead to non-negligible numerical error in orthogonality
\citep{saad2011numerical}. This was found to be the case especially for skewed
elements. The tolerance on the numerical error in orthogonality is chosen to be
$10^{-13}$ and if the error is larger than this value, an additional round of
reorthogonalization is performed. This rectified the issue and yielded columns
that are orthogonal up to machine precision.

\section{Remarks} \label{sec:remarks}

\subsection{Evaluating the divergence-free basis functions at a given set of
  points}

The divergence-free basis vectors
$\bm{\varphi}^{(e)}_1,\dots,\bm{\varphi}^{(e)}_n$ constructed in the previous
section contain the values of the constructed divergence-free basis functions at
the quadrature points in element $e$. The value of the basis functions at any
other point $\bm{s}^{(e)}$ in the simplex can be computed as follows. To develop
the algorithm, we temporarily revert back to the symbolic
notation. $\bm{\varphi}_{\ell}^{(e)}$ and $q_r$ now denote polynomials instead
of vectors. Each $\bm{\varphi}_{\ell}^{(e)}$ can be expanded in terms of the
polynomial basis function $q_r$ as
$\bm{\varphi}_{\ell}^{(e)}=\sum_{i=1}^d\sum_{r=1}^{p_j}N_{(i-1)p+r,\ell}^{(e)}q_r\left(\bm{x}\left(\bm{x}^{(e)}\right)\right)\bm{e}_i$. To
compute the value of $\bm{\varphi}_{\ell}^{(e)}$ at $\bm{s}^{(e)}$, we need the
value of $q_r$ at the mapped coordinate $\bm{s}(\bm{s}^{(e)})$. To compute this,
note that the polynomial $q_r$ is $1$ for $r=1$ and satisfies the below
Arnoldi-like relation for $r>1$:
\begin{equation*}
  x_iq_t=\sum_{r'=1}^{r}H_{r',r-1}q_{r'},
\end{equation*}
Here, $t=p_{j-1}-C_{d-i}^{j-1+d-i}+1,\dots,p_{j-1}$,
$r=p_{j-1}+\sum_{i'=1}^{i-1}C_{d-i'}^{j-1+d-i'}+1,\dots,p_{j-1}+\sum_{i'=1}^{i}C_{d-i'}^{j-1+d-i'}$,
$i=1,\dots,d$, and $j=1,\dots,k$. This relation also holds at $\bm{x}=\bm{s}$,
i.e., $s_iq_t(\bm{s})=\sum_{r'=1}^rH_{r',r-1}q_{r'}(\bm{s})$. Note that the
$H_{i,j}$s here are known quantities. Therefore, rearranging it yields the
expression for $q_r(\bm{s})$ to be
\begin{equation*}
  q_r(\bm{s})=\frac{1}{H_{r,r-1}}\left(s_iq_t-\sum_{r'=1}^{r-1}H_{r',r-1}q_{r'}\right)
\end{equation*}
for $r>1$. Note that this is a recursive expression. If we know the values of
$q_1(\bm{s}),\dots,q_r(\bm{s})$, then the value of $q_{r+1}(\bm{s})$ can be
computed using it. This yields the below recursive algorithm to compute
$q_r(\bm{s})$ for all $r$:
\begin{algorithmic}[1]
  \State $q_1(\bm{s})=1$
  \For{$j=1,\dots,k$}
    \State $c=p_{j-1}$
    \For{$i=1,\dots,d$}
      \State $j''=C_{d-i}^{j-1+d-i}$
      \For{$j'=1,\dots,j''$}
        \State $q_{c+1}(\bm{s})=\frac{1}{H_{c+1,c}}\left(s_iq_{p_{j-1}-j''+j'}(\bm{s})-\sum_{r'=1}^{c}H_{r',c}q_{r'}(\bm{s})\right)$
        \State $c=c+1$
      \EndFor
    \EndFor
  \EndFor
\end{algorithmic}
Using the obtained $q_r(\bm{s})$, the value of the divergence-free basis
functions can be computed by taking their linear combination. For a MATLAB
implementation of this algorithm, we refer the reader to figure
\ref{fig:ardivfreebfeval} in appendix \ref{app:matlabeval}.

\begin{figure}[t]
\centering
\begin{minipage}[c]{\linewidth}
\begin{subfigure}{0.25\linewidth}
\centering
\adjustbox{max width=\linewidth,trim=0.1cm 0.1cm 0.1cm 0.1cm,clip}{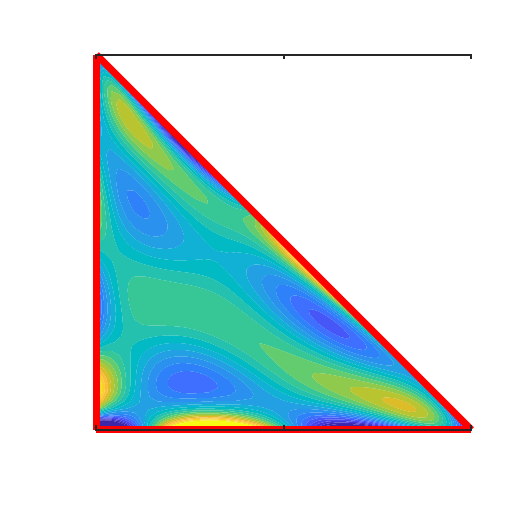}
\caption{}
\end{subfigure}\hfill
\begin{subfigure}{0.25\linewidth}
\centering
\adjustbox{max width=\linewidth,trim=0.1cm 0.1cm 0.1cm 0.1cm,clip}{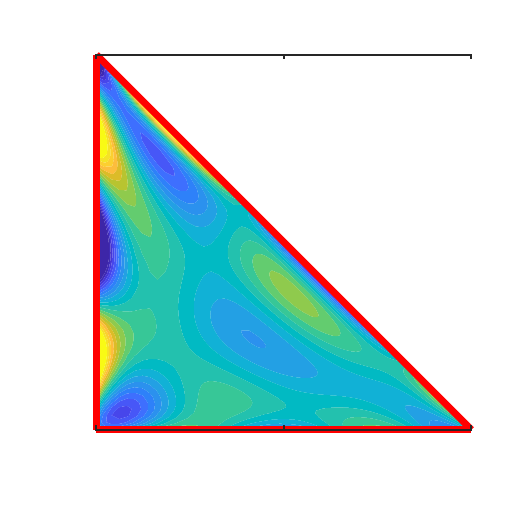}
\caption{}
\end{subfigure}\hfill
\begin{subfigure}{0.25\linewidth}
\centering
\adjustbox{max width=\linewidth,trim=0.1cm 0.1cm 0.1cm 0.1cm,clip}{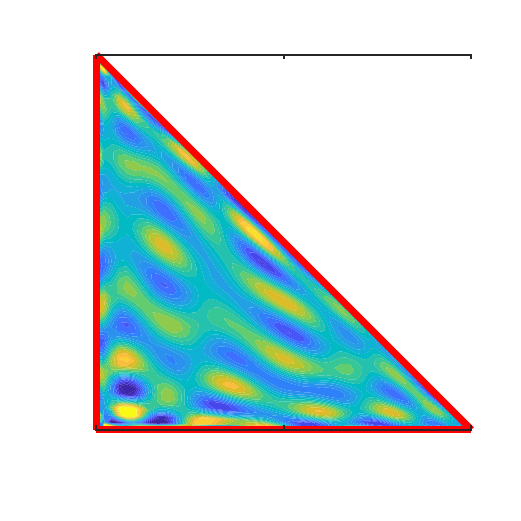}
\caption{}
\end{subfigure}\hfill
\begin{subfigure}{0.25\linewidth}
\centering
\adjustbox{max width=\linewidth,trim=0.1cm 0.1cm 0.1cm 0.1cm,clip}{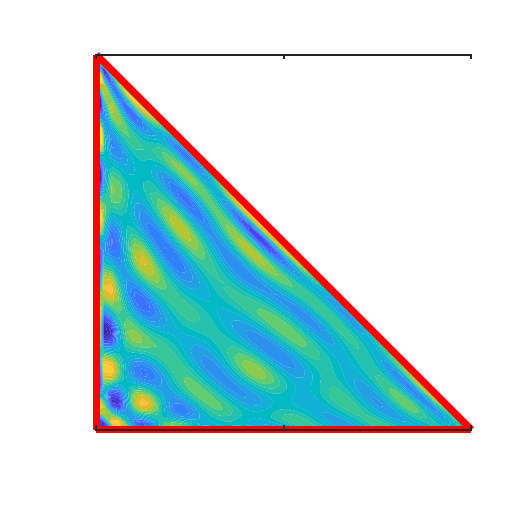}
\caption{}
\end{subfigure}\hfill
\begin{subfigure}{0.3\linewidth}
\centering
\adjustbox{max width=\linewidth,trim=0.2cm 0.3cm 0.9cm 0.5cm,clip}{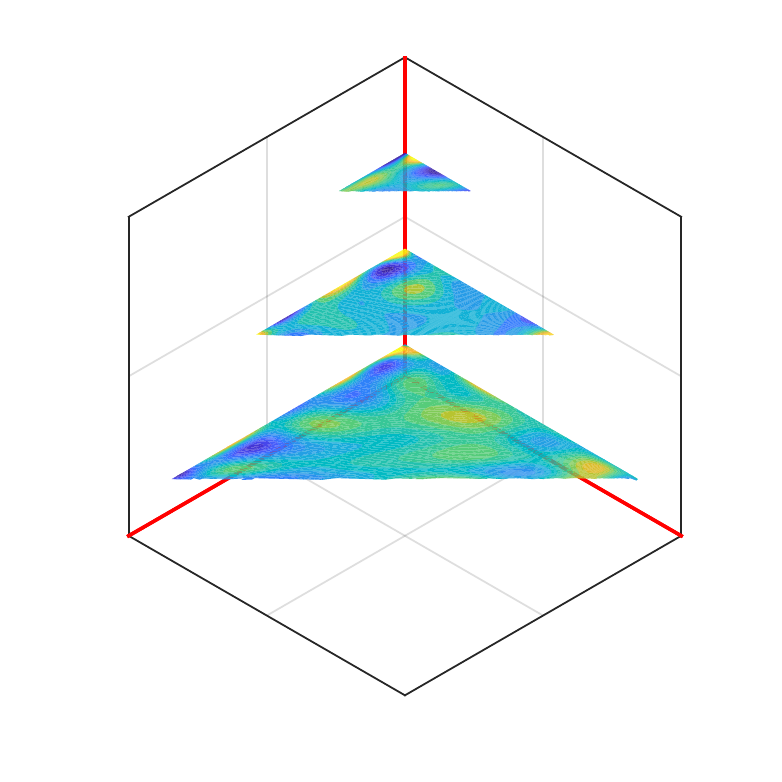}
\caption{}
\end{subfigure}
\begin{subfigure}{0.3\linewidth}
\centering
\adjustbox{max width=\linewidth,trim=0.2cm 0.3cm 0.9cm 0.5cm,clip}{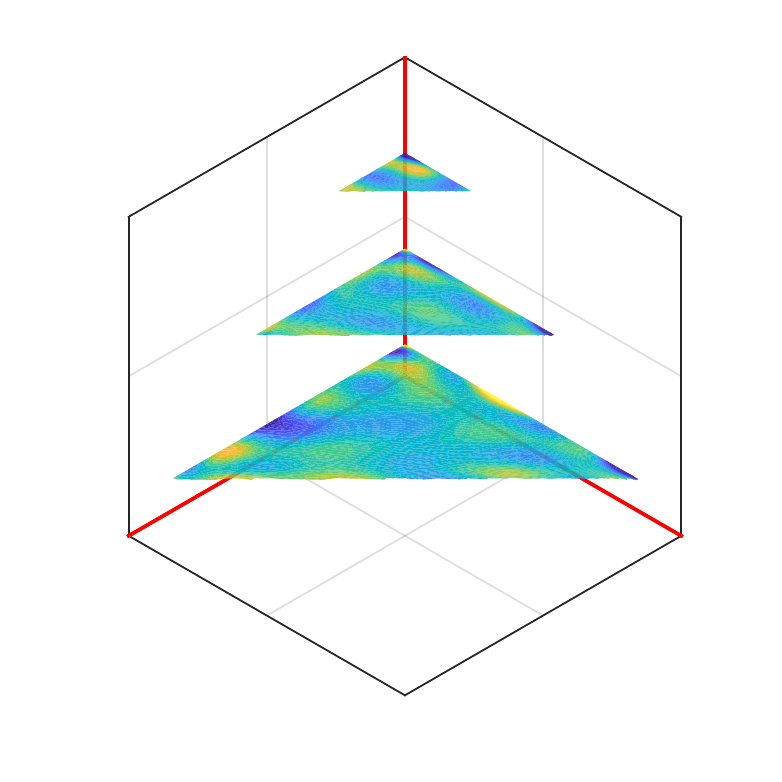}
\caption{}
\end{subfigure}
\begin{subfigure}{0.3\linewidth}
\centering
\adjustbox{max width=\linewidth,trim=0.2cm 0.3cm 0.9cm 0.5cm,clip}{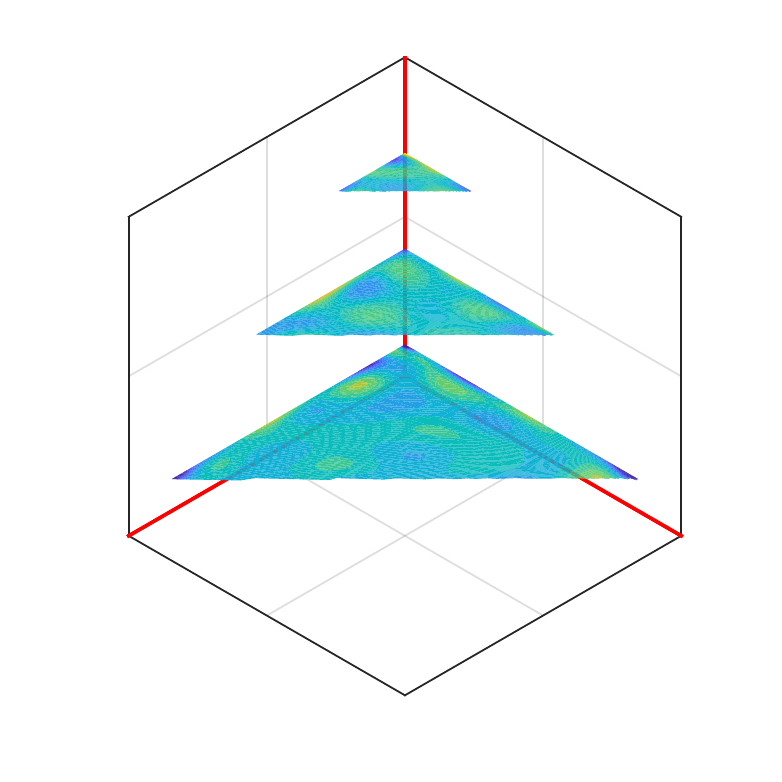}
\caption{}
\end{subfigure}
\centering
\begin{subfigure}{0.4\linewidth}
\centering
\adjustbox{max width=\linewidth,trim=0cm 6.2cm 0cm 0.5cm,clip}{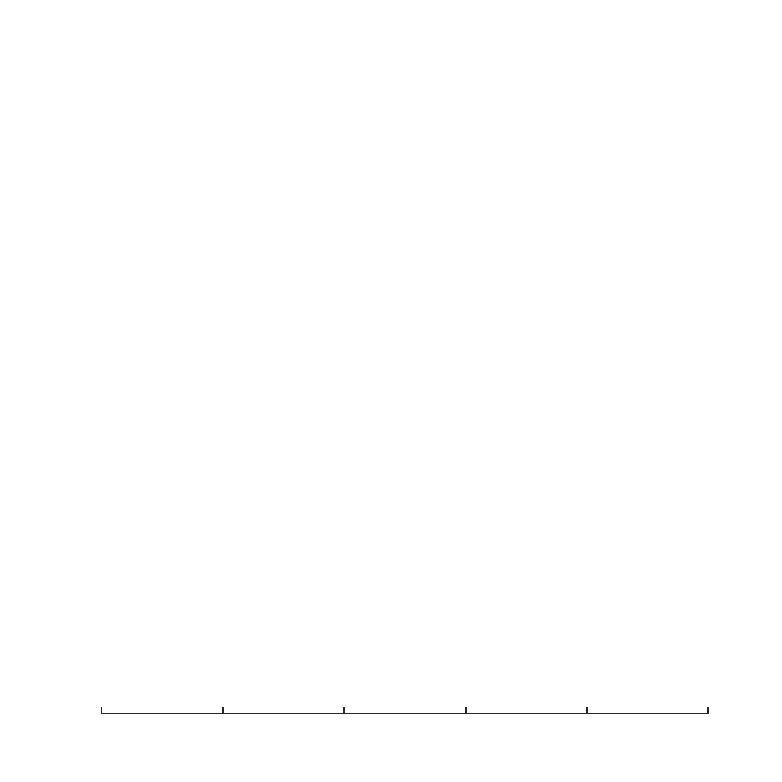}
\end{subfigure}
\end{minipage}
\caption{(a) $x$- and (b) $y$-component of a divergence-free basis function of degree $5$ in two dimensions. (c) $x$- and (d) $y$-component of a divergence-free basis function of degree $15$ in two dimensions. (e) $x$-, (f) $y$- and (g) $z$-component of a divergence-free basis function of degree 10 in three dimensions. The colorbar for all the contours is shown at the bottom.}
\label{fig:basis}
\end{figure}


\subsection{Plots of some divergence-free basis functions}

The contour plots of a few two- and three-dimensional divergence-free basis
functions constructed in the reference element are shown in figure
\ref{fig:basis}. Figures \ref{fig:basis}a and \ref{fig:basis}b show the $x$- and
$y$-components, respectively, of a two-dimensional basis function of degree
five. Figures \ref{fig:basis}c and \ref{fig:basis}d show the $x$- and
$y$-components, respectively, of another two-dimensional basis function of
degree 15. The $x$-, $y-$, and $z$-components of a three-dimensional basis
function of degree 10 are shown in figures \ref{fig:basis}e, \ref{fig:basis}f,
and \ref{fig:basis}g, respectively. Their slices at a few $z$-locations are
shown.


\begin{figure}[t]
\centering
\begin{subfigure}{0.3\linewidth}
\centering
\adjustbox{max width=\linewidth,trim=0cm 0cm 0cm 0cm,clip}{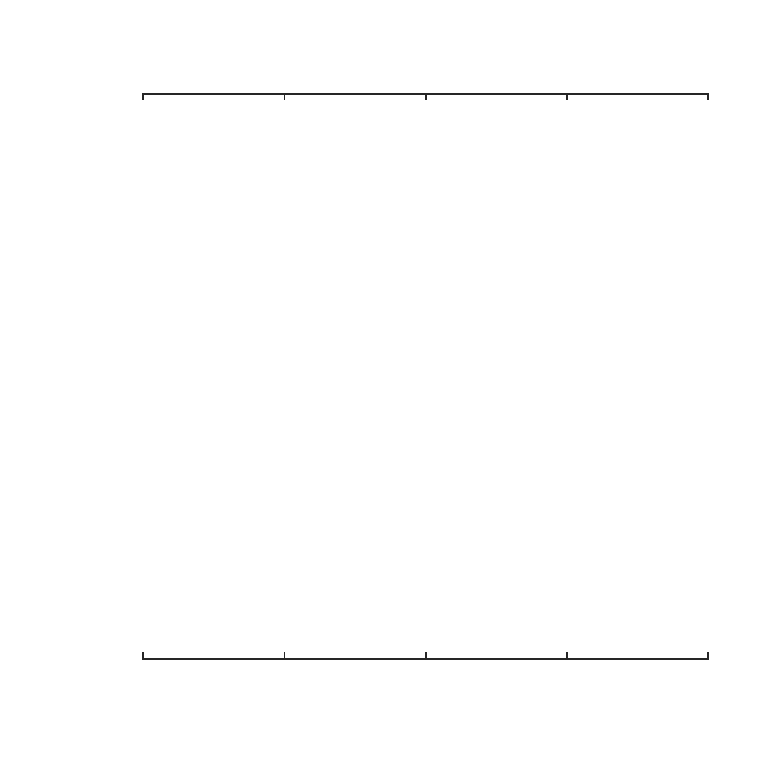}
\caption{}
\end{subfigure}
\begin{subfigure}{0.3\linewidth}
\centering
\adjustbox{max width=\linewidth,trim=0cm 0cm 0cm 0cm,clip}{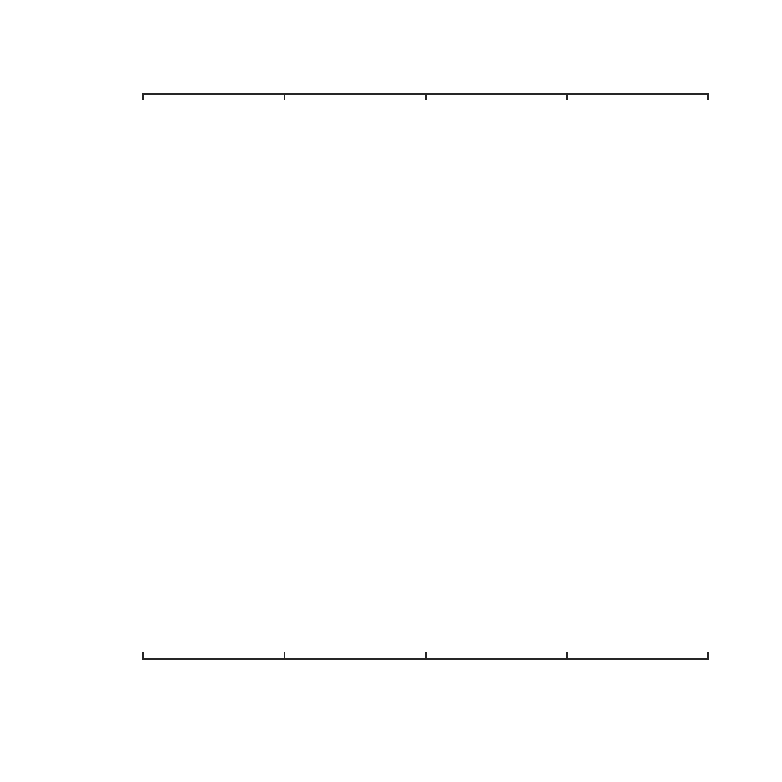}
\caption{}
\end{subfigure}
\begin{subfigure}{0.35\linewidth}
\centering
\adjustbox{max width=\linewidth,trim=0cm 0.7cm 0.7cm 1cm,clip}{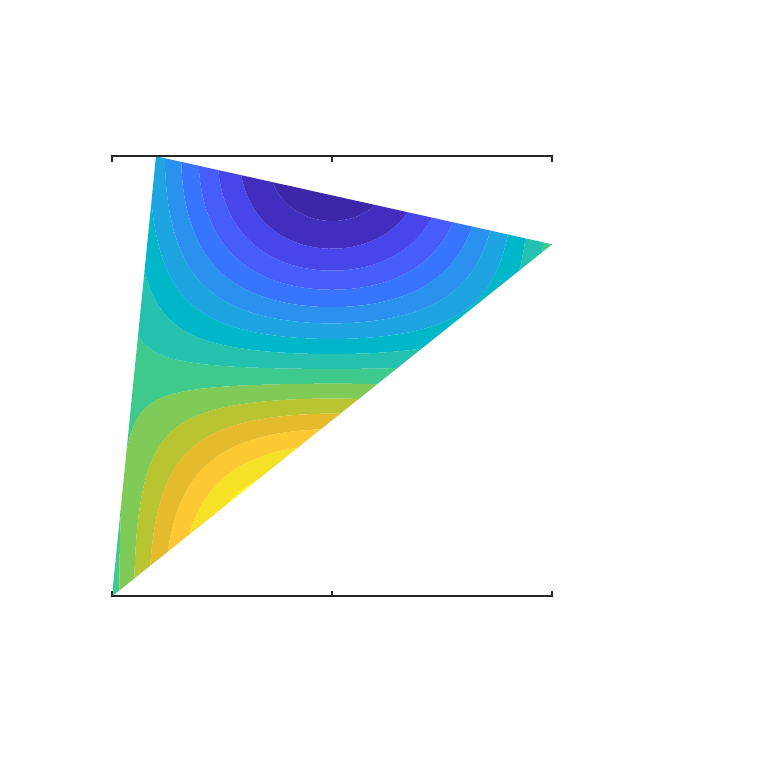}
\caption{}
\end{subfigure}
\caption{Divergence-free projection in a general triangle. (a) Projection error v/s polynomial degree ($\circ$ - $x$-component and $\times$ - $y$-component). (b) Constraint error v/s polynomial degree. (c) Contours of the $x$-component of the projected velocity field for polynomial degree $20$.}
\label{fig:divfree_bases_one_elem}
\end{figure}

\subsection{Example divergence-free projection in a general triangle}

In this example, we project the divergence-free velocity field
$\bm{u}(x,y)=(\sin(\pi x)\cos(\pi y),-\cos(\pi x)\sin(\pi y))$ (Taylor-Green
velocity field) onto the divergence-free basis constructed in a general
triangle.  The node-coordinate matrix of the triangle is
$X^{(e)}=[0\hspace{.5em}0;1\hspace{.5em}0.8;0\hspace{.5em}0.1]$. The MATLAB code
used to perform the projection is given below:
\begin{small}
\begin{verbatim}
k=20;d=2;Xe=[0 0;1 0.8;0.1 1];
tgfac=pi;ffu=@(x,y)sin(tgfac*x).*cos(tgfac*y);ffv=@(x,y)-cos(tgfac*x).*sin(tgfac*y);
[N,Q,H,Qd,x,w,~,C]=ardivfreebfref(k,d);
[Ne,Qde]=ardivfreebfgen(k,d,Xe,N,Q);
kp1d=(k+1)^d;xe=repmat(Xe(1,:),size(x(:,1)));xe=xe+x*(Xe(2:d+1,:)-repmat(Xe(1,:),[d 1]));
f=[ffu(xe(:,1),xe(:,2)); ffv(xe(:,1),xe(:,2))];
wd=repmat(w,[d 1]);wf=factorial(d);
dc=mgs_with_reorth(Qde,f,wd,wf,size(Qde,2),1);
\end{verbatim}
\end{small}
Here, \texttt{mgs\_with\_reorth} is the function in figure
\ref{fig:ardivfreebfref} and \texttt{dc} is the vector storing the coefficients
of the projection. The error in projection and in the satisfaction of the
divergence-free constraint are computed using the below MATLAB code.
\begin{small}
\begin{verbatim}
nplt=50;npltd=nplt^d;s1{1}=linspace(0,1,nplt); for j=2:d, s1{j}=s1{1}; end;
[st{1:d}]=ndgrid(s1{:}); s=reshape(cat(d,st{:}),[npltd d]);
s=[s(:,1) s(:,2:d).*cumprod(1-s(:,1:d-1),2)];
[Wde,W]=ardivfreebfeval(k,d,H,Ne,s);
se=repmat(Xe(1,:),size(s(:,1)));se=se+s*(Xe(2:d+1,:)-repmat(Xe(1,:),[d 1]));
fex=[ffu(se(:,1),se(:,2)) ffv(se(:,1),se(:,2))];
for i=1:d, err{i}=[]; end; yy=zeros(d*npltd,1);jdimp=0;jddimp=0;
errdiv=[0];kdim=nchoosek(k+d,d);
Xet=Xe'; Fe=Xet(:,2:d+1)-repmat(Xet(:,1),[1 d]); clear Xet;
Feinv=kron(inv(Fe),speye(nchoosek(k+d,d))); Ce=C*Feinv;
for j=0:k
  jdim=nchoosek(j+d,d); jddim=d*jdim-jdimp;
  yy=yy+Wde(:,jddimp+1:jddim)*dc(jddimp+1:jddim);
  for i=1:d
    err{i}=[err{i} max(abs(yy((i-1)*npltd+1:i*npltd)-fex(:,i)))];
  end
  ii=[]; for i=1:d, ii=[ii (i-1)*kdim+1:(i-1)*kdim+jdim]; end      
  errdiv=[errdiv max(abs(Ce(1:jdimp,ii)*Ne(ii,1:jddim)*dc(1:jddim)))];
  jdimp=jdim;jddimp=jddim;
end
\end{verbatim}
\end{small}
The projection error is defined to be the maximum absolute difference between
$\bm{u}$ and the projection $\bm{u}_h$ evaluated at 2500 points in the
triangle. The constraint error is defined for polynomial degree $j$ as the
maximum absolute projection of the divergence of $\bm{u}_h$ along the
orthonormal polynomials of degree less than $j$ (
$\max_{1\leq j'\leq
  C_d^{j-1+d}}|\int_{\Omega^{(e)}}q_{j'}\bm{\nabla}\cdot\bm{u}_h\,d\Omega|$).

Figures \ref{fig:divfree_bases_one_elem}a and \ref{fig:divfree_bases_one_elem}b
show the computed error in projection and in the satisfaction of the integral
divergence-free constraint, respectively, as a function of polynomial
degree. The projection error decreases exponentially with increasing polynomial
degree and reaches machine precision for degree 19. The integral divergence-free
constraint is satisfied up to machine precision for all polynomial degrees. The
contours of the $x$-component of the projected vector computed with $k=20$ are
shown in figure \ref{fig:divfree_bases_one_elem}c.

Note that the coefficients of projection (stored in the vector \texttt{dc}) are
computed only for polynomial degree 20. Since the constructed basis functions
are hierarchical, the coefficients for degree $j$ less than 20 are nothing but
the first $dC_d^{j+d}-C_d^{j-1+d}$ entries of the vector \texttt{dc}. Therefore,
the projected function \texttt{yy} is incrementally computed at the points
\texttt{s} for each polynomial degree \texttt{j} as {\small
  \texttt{yy=yy+Wde(:,jddimp+1:jddim)*dc(jddimp+1:jddim)}}.

\subsection{Computational cost}

The computational cost to generate the divergence-free basis functions in the
reference element (step 1) scales as $O((k+1)^{3d})$. The cost to construct the
basis functions function in a general element (step 2) also scales as
$O((k+1)^{3d})$, thought with a much smaller constant. The cost of evaluating
the basis functions at $n_p$ points with \texttt{ardivfreebfeval} scales as
$O(n_p(k+1)^{2d})$.

The most expensive part in the construction is step 1. However, we note that step 1 needs to be performed just once for the largest polynomial degree of interest, say $k$. The basis functions for degrees smaller than $k$ are part of the degree $k$ basis because of the hierarchy of the basis functions. The outputs of step 1 can even be precomputed and stored in a file that can be read at the beginning of each simulation. We note that step 2, which needs to be performed for each element in the mesh, is substantially cheaper than step 1. For example, for polynomial degree 15 in three dimensions, step 1 takes 20 seconds, while step 2 consumes just 0.3 seconds. Therefore, our methodology is a computationally efficient procedure to compute an orthonormal and hierarchical divergence-free basis for multiple elements.

\subsection{On the structure of matrices $N$, $N^{(e)}$, and $H$}

\begin{figure}[t]
\centering
\begin{subfigure}{0.2\linewidth}
\centering
\adjustbox{max width=\linewidth,trim=1.2cm 0cm 1.5cm 0cm,clip}{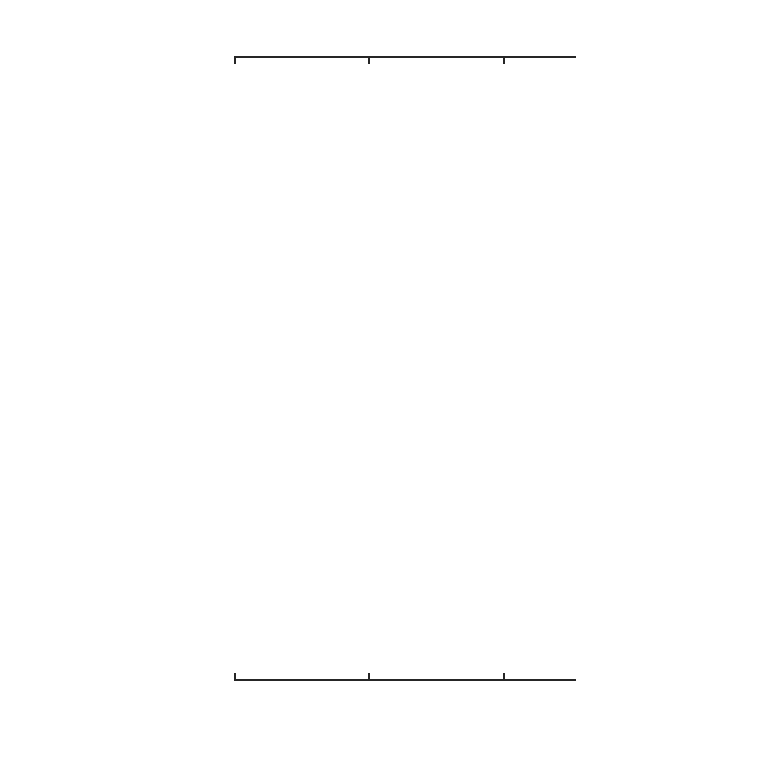}
\caption{}
\end{subfigure}
\begin{subfigure}{0.3\linewidth}
\centering
\adjustbox{max width=\linewidth,trim=0cm 0cm 0cm 0cm,clip}{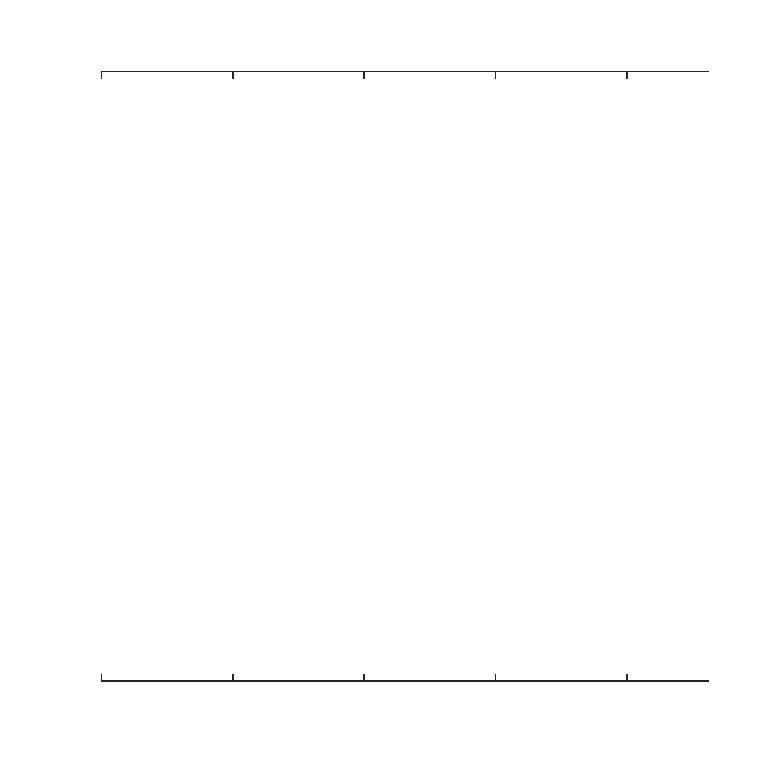}
\caption{}
\end{subfigure}
\caption{(a) Sparsity of $N$ and $N^{(e)}$. (b) Sparsity of $H$ with entries smaller than $10^{-13}$ neglected.}
\label{fig:sparsity}
\end{figure}

Nearly half of the entries in the coefficient matrices $N$ and $N^{(e)}$ are
zeros and their non-zero (sparsity) patterns are identical. Consider the
expression for the divergence-free basis functions $\bm{\varphi}_{\ell}$ and
$\bm{\varphi}_{\ell}^{(e)}$ of degree $j$:
$\bm{\varphi}_{\ell}=\sum_{i=1}^d\sum_{r=1}^{p_j}N_{(i-1)p+r,\ell}q_r\bm{e}_i$
and
$\bm{\varphi}_{\ell}^{(e)}=\sum_{i=1}^d\sum_{r=1}^{p_j}N_{(i-1)p+r,\ell}^{(e)}q_r\left(\bm{x}\left(\bm{x}^{(e)}\right)\right)\bm{e}_i$,
where $\ell=n_{j-1}+1,\dots,n_j$. Notice that the summation along $r$ is from
$1$ to $p_j$ and not from $1$ to $p$, and therefore, the corresponding entries
in matrices $N$ and $N^{(e)}$ are zeros. Specifically, $N_{(i-1)p+r,\ell}=0$ and
$N_{(i-1)p+r,\ell}^{(e)}=0$ for $r$ larger than $p_j$ but less than or equal to
$p$ and for each $i=1,\dots,d$. This non-zero pattern of $N$ and $N^{(e)}$ are
shown in figure \ref{fig:sparsity}a. The matrices computed for the
divergence-free projection problem discussed previously are used to plot this
figure. Notice the staircase pattern of the matrices. This non-zero pattern can
be used to reduce the cost of multiplying vectors with the matrices $N$ and
$N^{(e)}$.




The structure of the upper-Hessenberg matrix $H$ generated using the
Arnoldi-based procedure is very interesting. Several of its entries in the upper
triangular portion are very close to zero. To show this, we define $H_z$ to be
the matrix which is same as $H$ except that the entries in $H$ that are smaller
than $10^{-13}$ are set to zero in $H_z$. Figure \ref{fig:sparsity}b shows the
sparsity of $H_z$ for polynomial degree 20 in two dimensions. The sparsity
pattern has three diverging bands comprising of block matrices. This is not
accidental. It is a multi-dimensional analogue of the three-term recurrence
relation of the one-dimensional orthonormal polynomials. One might be tempted to
exploit this pattern to develop an Arnoldi-based process that would cost
$O(k+1)^{2d}$ number of operations instead of the current $O(k+1)^{3d}$ cost. We
developed one such method. But the algorithm was numerically unstable. The
generated polynomials lost orthogonality. This is because of the same reason the
Lanczos vectors lose orthogonality in finite-precision arithmetic without
selective or complete reorthogonalization \citep{saad2011numerical}. Note that
in the Arnoldi-based process in \textbf{step 1}, each vector $q_j$ is made
orthogonal to all previous vectors $q_1,\dots,q_{j-1}$. This is essentially
complete reorthogonalization and it is necessary to retain orthogonality of the
vectors $q_j$.

\section{Applications} \label{sec:appln}

We use the constructed divergence-free basis functions to compute numerical
solutions of some PDEs. For the first application, we show in detail how to
exploit the orthonormal and hierarchical features of our basis. An efficient
implementation of the hybridized mixed method is presented to compute numerical
solutions for all polynomial degrees from zero to a given $k$. For the remaining
applications, just the results are presented. Efficient implementations can be
constructed similarly.

\subsection{Helmholtz projection}

The problem is as follows: Consider a domain $\Omega$. Let $T_h$ denote a
triangulation of this domain. Given a function $\bm{g}$, compute its projection
onto the divergence-free basis constructed in each element of the triangulation
such that the normal component of the projection is continuous across the
inter-element boundaries. This problems amounts to solving the below
Helmholtz-type PDE problem:
\begin{alignat}{3}
&\bm{u}\quad+&\bm{\nabla}\lambda&=\bm{g}&&\quad\text{ in }\Omega,\\
&\bm{\nabla}\cdot \bm{u}& &=0&&\quad\text{ in }\Omega,\text{ and}\nonumber\\
&&\lambda&=0&&\quad\text{ on }\partial\Omega.\nonumber
\end{alignat}
Here, $\bm{u}$ is the desired projection, $\partial \Omega$ is the boundary of
$\Omega$, and $\lambda$ is the Lagrange multiplier that imposes the
divergence-free condition on $\bm{u}$. These equations are solved using the
hybridized BDM mixed method \citep{brezzi1985two} with one modification. The
proposed divergence-free basis is used in place of the usual polynomial basis to
approximate $\bm{u}$ in each element.




In the usual hyridized BDM mixed method for the above problem, $\bm{u}$ and $\lambda$ are approximated in each element $e$ by a polynomial of degree less than or equal to $k$ and $k-1$, respectively. Denote these approximations by $\bm{u}_h^{(e)}$ and $\lambda_h^{(e)}$. On each interior face $f$ of the mesh, the continuity of the normal component of the discontinuous approximation $\bm{u}_{h}^{(e)}$ is enforced using a Lagrange multiplier $\what{\lambda}^{(f)}_h$.  This Lagrange multiplier is taken to be a polynomial of degree less than or equal $k$ on each face $f$ and is an approximation to $\lambda$ on the faces. In each element $e$, $\bm{u}_h^{(e)}$ and $\lambda_h^{(e)}$ are defined to be solution to the problem:
\begin{alignat}{3} \label{eqn:ulambdae}
\ipiv{\bm{u}_{h}^{(e)}}{\bm{v}}{\Omega^{(e)}}&-\ipi{\lambda_{h}^{(e)}}{\bm{\nabla}\cdot\bm{v}}{\Omega^{(e)}}&&=\ipiv{\bm{g}}{\bm{v}}{\Omega^{(e)}}-\sum_{f\in F(e)}\ipibdom{\what{\lambda}_h^{(f)}}{\bm{v}\cdot\bm{n}}{\Gamma^{(f)}} &&\quad \forall \bm{v} \in [P_k(\Omega^{(e)})]^d,\\
\ipi{w}{\bm{\nabla}\cdot\bm{u}_{h}^{(e)}}{\Omega^{(e)}}&\,&&=0 &&\quad \forall w\in P_{k-1}(\Omega^{(e)}),\nonumber
\end{alignat}
for each element $e$. Here, $\bm{v}$ and $w$ are test functions, $F(e)$ is the set of faces of element $e$, and $\Gamma^{(f)}$ is the domain of face $f$. The equations for the Lagrange multiplier $\what{\lambda}_h^{(f)}$ are the normal continuity constraints:
\begin{equation} \label{eqn:gloprob1}
\int_{\Gamma^{(f)}}\left(\bm{u}_h^{(e^+)}\cdot \bm{n}^+ +\bm{u}_h^{(e^-)}\cdot \bm{n}^-\right)\mu\,d\Gamma=0\quad \forall\mu\in P_k(\Gamma^{(f)}),
\end{equation}
for each interior face $f$ of the mesh. Here, $e^{+}$ and $e^{-}$ are elements adjacent to face $f$. $\bm{n}^+$ and $\bm{n}^-$ are unit vectors normal to face $f$ and outward to the elements $e^{+}$ and $e^{-}$, respectively. On the boundary faces $f$, $\what{\lambda}_h^{(f)}$ is set to zero in accordance with the boundary condition.



Note that $\bm{u}_h^{(e)}$ is divergence-free at each point in the element $e$
because $\bm{\nabla}\cdot[P_k
(\Omega^{(e)})]^d=P_{k-1}(\Omega^{(e)})$. Therefore, instead of approximating
$\bm{u}_h^{(e)}$ with polynomials of degree less than or equal to $k$, it can be
approximated with divergence-free polynomials of degree less than or equal to
$k$. This simplifies \citep{cockburn2009two} equation \ref{eqn:ulambdae} to the
below projection problem:
\begin{alignat}{4}
\ipiv{\bm{u}_{h}^{(e)}}{\bm{v}}{\Omega^{(e)}}&&&=\ipiv{\bm{g}}{\bm{v}}{\Omega^{(e)}}-\sum_{f\in F(e)}\ipibdom{\what{\lambda}_h^{(f)}}{\bm{v}\cdot\bm{n}}{\Gamma^{(f)}} &&\quad \forall \bm{v} \in \bm{V}_k(\Omega^{(e)}),\label{eqn:simplocprob2_1}
\end{alignat}
where $\bm{V}_k(\Omega^{(e)})$ is the set of all divergence-free polynomials of
degree less than or equal to $k$ in element $e$. Note that $\lambda_h^{(e)}$ has
disappeared from the equation because $\bm{\nabla}\cdot\bm{v}=0$ for all test
functions $\bm{v}$ in $\bm{V}_k(\Omega^{(e)})$. To obtain an equation for just
$\what{\lambda}^{(f)}$, following \cite{cockburn2016static}, we decompose
$\bm{u}_h^{(e)}$ as
$\bm{u}^{(e)}_{h}=\bm{u}^{(e)}_{\bm{g}}+\bm{u}^{(e)}_{\what{\lambda}_h}$, where
$\bm{u}^{(e)}_{\bm{g}}$ and $\bm{u}^{(e)}_{\what{\lambda}_h}$ are defined as
solutions to the problems:
\begin{alignat}{4}
\ipiv{\bm{u}_{\bm{g}}^{(e)}}{\bm{v}}{\Omega^{(e)}}&=\ipiv{\bm{g}}{\bm{v}}{\Omega^{(e)}}&&\quad \forall\bm{v}\in\bm{V}_{k}(\Omega^{(e)}),\text{ and} \label{eqn:simplocprob1}
\end{alignat}
\begin{alignat}{4}
\bm{u}_{\what{\lambda}_h}^{(e)}=\sum_{f\in F(e)}\bm{u}_{\what{\lambda}^{(f)}_h}^{(e)},\text{ where }\ipiv{\bm{u}_{\what{\lambda}_h^{(f)}}^{(e)}}{\bm{v}}{\Omega^{(e)}}&=-\ipibdom{\what{\lambda}_h^{(f)}}{\bm{v}\cdot\bm{n}}{\Gamma^{(f)}}&&\quad \forall \bm{v}\in\bm{V}_{k}(\Omega^{(e)}),\label{eqn:simplocprob2}
\end{alignat}
respectively. Then, we substitute them into equation \ref{eqn:gloprob1} and this yields the below desired equation for just $\what{\lambda}_h^{(f)}$:
\begin{equation} \label{eqn:gloprob}
\begin{split}
\sum_{e\in\{e^+,e^-\}}\sum_{f'\in F(e)}\ipiv{\bm{u}^{(e)}_{\what{\lambda}^{(f')}_h}}{\bm{u}^{(e)}_{\mu}}{\Omega^{(e)}}=-\sum_{e\in\{e^+,e^-\}}\ipiv{\bm{u}^{(e)}_{\bm{g}}}{\bm{u}^{(e)}_{\mu}}{\Omega^{(e)}}\quad
\forall\mu\in P_k(\Gamma^{(f)}).
\end{split}
\end{equation}



To obtain the corresponding matrix problems, we expand $\bm{u}^{(e)}_{\bm{g}}$ and $\bm{u}^{(e)}_{\what{\lambda}_h^{(f)}}$ in terms of the divergence-free basis functions constructed in each element $e$ and expand $\what{\lambda}_h^{(f)}$ on each face $f$ in terms of  the Arnoldi-based orthonormal polynomials constructed in the reference face element:
\begin{equation}\label{eqn:expand}
\what{\lambda}_h^{(f)}=\sum_{j=1}^{\wtilde{m}}c_j^{(f)}q_j^{(f)},\,\bm{u}_{\bm{g}}^{(e)}=\sum_{\ell=1}^{n}\alpha^{(e)}_\ell\bm{\varphi}^{(e)}_{\ell},\,\bm{u}_{\what{\lambda}_h^{(f)}}^{(e)}=\sum_{j=1}^{\wtilde{m}}c_j^{(f)}\bm{u}_{q_j^{(f)}}^{(e)},
\end{equation}
where $\wtilde{m}$ is the dimension of $P_k(\Gamma^{(f)})$ (which equals $C_{d-1}^{k+d-1}$), and $\bm{u}_{q_j^{(f)}}^{(e)}=\sum_{\ell=1}^{n}\beta^{(e,f)}_{\ell,j}\bm{\varphi}^{(e)}_{\ell}$ is the solution to equation \ref{eqn:simplocprob2} with $\what{\lambda}_h^{(f)}$ set to $q_j^{(f)}$, i.e.,
\begin{alignat}{4}
\ipiv{\bm{u}_{q_j^{(f)}}^{(e)}}{\bm{v}}{\Omega^{(e)}}&=-\ipibdom{q_j^{(f)}}{\bm{v}\cdot\bm{n}}{\Gamma^{(f)}}&&\quad \forall \bm{v}\in\bm{V}_{k}(\Omega^{(e)}),\label{eqn:simplocprobqj}
\end{alignat}
and $c_j^{(f)}$, $\alpha^{(e)}_\ell$ and $\beta^{(e,f)}_{\ell,j}$ are the coefficients. Substituting the expressions for $\bm{u}_{\bm{g}}^{(e)}$ and $\bm{u}_{q_j^{(f)}}^{(e)}$ in equations \ref{eqn:simplocprob1} and \ref{eqn:simplocprobqj}, respectively, and requiring the equality for the test function $\bm{v}$ equal to each divergence-free basis function results in the below matrix problems for the coefficients $\alpha^{(e)}_\ell$ and $\beta^{(e,f)}_{\ell,j}$:
\begin{alignat}{4}
\sum_{\ell=1}^n\left(\ipiv{\bm{\varphi}^{(e)}_\ell}{\bm{\varphi}^{(e)}_i}{\Omega^{(e)}}\right)\alpha_{\ell}^{(e)}&=\ipiv{\bm{g}}{\bm{\varphi}_i^{(e)}}{\Omega^{(e)}}&&\quad\text{for }i=1,\dots,n,\text{ and} \label{eqn:simplocprob1mat}\\
\sum_{\ell=1}^n\left(\ipiv{\bm{\varphi}^{(e)}_{\ell}}{\bm{\varphi}^{(e)}_i}{\Omega^{(e)}}\right)\beta^{(e,f)}_{\ell,j}&=-\ipibdom{q_{j}^{(f)}}{\bm{\varphi}^{(e)}_i\cdot\bm{n}}{\Gamma^{(f)}}&&\quad\text{for }i=1,\dots,n;j=1,\dots,\wtilde{m};\text{ and }\forall f\in F(e).\label{eqn:simplocprob2mat}
\end{alignat}
The above two problems are called the local problems. They need to be solved in
each element of the mesh. Substituting the expression for
$\bm{u}_{\what{\lambda}_h^{(f)}}^{(e)}$ in equation \ref{eqn:gloprob} and
requiring the equality for $\mu$ equal to each orthonormal polynomial yields the
below equations for the coefficients $c_j^{(f)}$
\begin{equation} \label{eqn:glomatprob}
\begin{split}
\sum_{e\in\{e^+,e^-\}}\sum_{f'\in F(e)}\sum_{j=1}^{\wtilde{m}}\left(\ipiv{\bm{u}^{(e)}_{q_j^{(f')}}}{\bm{u}^{(e)}_{q_{i}^{(f)}}}{\Omega^{(e)}}\right)c_j^{(f')}=-\sum_{e\in\{e^+,e^-\}}\ipiv{\bm{u}^{(e)}_{\bm{g}}}{\bm{u}^{(e)}_{q_i^{(f)}}}{\Omega^{(e)}},
\end{split}
\end{equation}
for $i=1,\dots,\wtilde{m}$ and for each interior face $f$ of the mesh. The above
problem is called the global problem. It couples the individual local problems.
 
The orthonormality and hierarchial features of the basis functions simplfy the
local problem solution and the global problem assembly. Since the proposed basis
functions $\bm{\varphi}_{\ell}^{(e)}$ are orthonormal, the local problems
simplify to computing the below inner-products:
\begin{alignat}{4}
\alpha_{\ell}^{(e)}&=\left(\ipiv{\bm{g}}{\bm{\varphi}_{\ell}^{(e)}}{\Omega^{(e)}}\right)/|{\Omega^{(e)}}|&&\quad\text{for }\ell=1,\dots,n,\text{ and} \label{eqn:simplocprob1mat1}\\
\beta^{(e,f)}_{\ell,j}&=-\left(\ipibdom{q_{j}^{(f)}}{\bm{\varphi}^{(e)}_{\ell}\cdot\bm{n}}{\Gamma^{(f)}}\right)/|{\Omega^{(e)}}|&&\quad\text{for }\ell=1,\dots,n;j=1,\dots,\wtilde{m};\text{ and }\forall f\in F(e).\label{eqn:simplocprob2mat1}
\end{alignat}
Substituting the expression for $\bm{u}^{(e)}_{q_{i}^{(f)}}$ and using the
orthogonality of the divergence-free basis functions
$\bm{\varphi}_{\ell}^{(e)}$s simplifies the global problem to:
\begin{equation} \label{eqn:glomatprob1}
\begin{split}
\sum_{e\in\{e^+,e^-\}}\sum_{f'\in F(e)}\sum_{j=1}^{\wtilde{m}}\left(\sum_{\ell=1}^n\beta^{(e,f')}_{\ell,j}\beta^{(e,f)}_{\ell,i}|\Omega^{(e)}|\right)c_j^{(f')}=-\sum_{e\in\{e^+,e^-\}}\sum_{\ell=1}^n\alpha_{\ell}^e\beta^{(e,f)}_{\ell,i}|\Omega^{(e)}|,
\end{split}
\end{equation}
for $i=1,\dots,\wtilde{m}$ and for each interior face $f$ of the mesh. These
equations can be written as the matrix problem $Ax=b$. Here, $A$ is the
left-hand side matrix of size $n_f\wtilde{m}\times n_f\wtilde{m}$, where $n_f$
is the number of interior faces in the mesh. $x$ is the vector of coefficients
of size $n_f\wtilde{m}\times 1$ and is defined as
$x_{(f-1)\wtilde{m}+i}=c_i^{(f)}$. $b$ is the right-hand side vector of size
$n_f\wtilde{m}\times 1$. The entries of the left-hand side matrix and the
right-hand side vector are defined as:
\begin{equation*}
A_{(f-1)\wtilde{m}+i,(f'-1)\wtilde{m}+j}=\sum_{e\in\{e^+,e^-\}}\sum_{\ell=1}^n\beta^{(e,f')}_{\ell,j}\beta^{(e,f)}_{\ell,i}|\Omega^{(e)}|\text{ and }b_{(f-1)\wtilde{m}+i}=-\sum_{e\in\{e^+,e^-\}}\sum_{\ell=1}^n\alpha_{\ell}^{(e)}\beta^{(e,f)}_{\ell,i}|\Omega^{(e)}|.
\end{equation*}
Note that the left-hand side matrix is sparse \citep{cockburn2016static} because $A_{(f-1)\wtilde{m}+i,(f'-1)\wtilde{m}+j}$ is non-zero if and only if the faces $f$ and $f'$ belong to a common element.



The hierarchical feature of the basis functions can be exploited to develop an
efficient assembly procedure for all polynomial degree from zero to $k$. The
left-hand side matrix and the right-hand side vector are assembled using element
matrix $A^{(e)}$ (size $(d+1)\wtilde{m}\times(d+1)\wtilde{m}$) and element
vector $b^{(e)}$ (size $(d+1)\wtilde{m}\times 1$) that are computed for each
element $e$ as:
\begin{equation}
\begin{split}
A^{(e)}_{(g-1)\wtilde{m}+i,(g'-1)\wtilde{m}+j}=\sum_{\ell=1}^n\beta^{(e,g')}_{\ell,j}\beta^{(e,g)}_{\ell,i}|\Omega^{(e)}|\text{ and }b^{(e)}_{(g-1)\wtilde{m}+i}=\sum_{\ell=1}^n\alpha_{\ell}^{(e)}\beta^{(e,g)}_{\ell,i}|\Omega^{(e)}|.
\end{split}
\end{equation}
Here, $g$ and $g'$ are the local index of the $(d+1)$ faces of the
element. Since the basis functions are hierarchical, for each polynomial degree
$k'$ less than $k$, the element matrix $A^{(e,k')}$ (size
$(d+1)\wtilde{m}'\times(d+1)\wtilde{m}'$) and the element vector $b^{(e,k')}$
(size $(d+1)\wtilde{m}'\times 1$) are the following partial sums of the
summations in the above equation:
\begin{equation}
\begin{split}
A^{(e,k')}_{(g-1)\wtilde{m}'+i,(g'-1)\wtilde{m}'+j}=\sum_{\ell=1}^{n'}\beta^{(e,g')}_{\ell,j}\beta^{(e,g)}_{\ell,i}|\Omega^{(e)}|\text{ and }b^{(e)}_{(g-1)\wtilde{m}'+i}=\sum_{\ell=1}^{n'}\alpha_{\ell}^{(e)}\beta^{(e,g)}_{\ell,i}|\Omega^{(e)}|.
\end{split}
\end{equation}
where $n'=dC_d^{k'+d}-C_d^{k'-1+d}$ and $\wtilde{m}'=C_{d-1}^{k'+d-1}$. Hence,
given the coefficients $\beta_{\ell,j}^{(e,g)}$ and $\alpha_{\ell}^{(e)}$ that
are computed for degree $k$, the element matrices and vectors for all degrees
$k'$ (up to $k$) can be incrementally constructed by performing an update to the
element matrix and vector computed for degree $k'-1$. Therefore, the left-hand
side matrix $A^{(k')}$ and the right-hand side vector $b^{(k')}$ for all degrees
$k'$ (up to $k$) can also be incrementally assembled by updating the matrix
$A^{(k'-1)}$ and the vector $b^{(k'-1)}$ computed for degree $k'-1$. Using this
idea, all the left-hand side matrices and the right-hand side vectors from
polynomial degree $0$ to the given degree $k$ are efficiently
constructed. 

\begin{figure}[t]
\centering
\begin{subfigure}[b]{0.35\linewidth}
\centering
\adjustbox{max width=\linewidth,trim=0cm 0.7cm 0.7cm 1cm,clip}{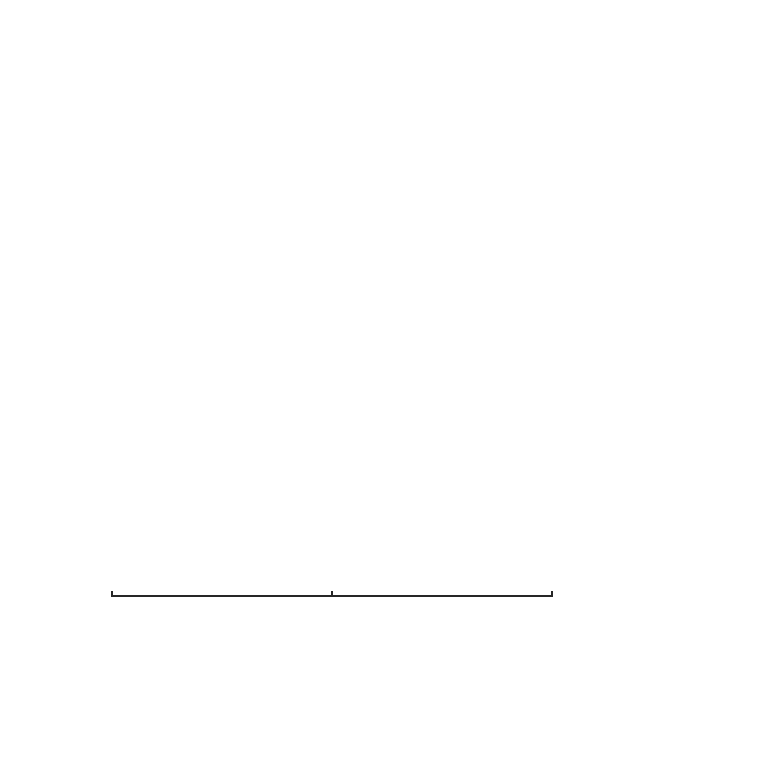}
\caption{}
\end{subfigure}
\begin{subfigure}[b]{0.3\linewidth}
\centering
\adjustbox{max width=\linewidth,trim=0cm 0cm 0cm 0cm,clip}{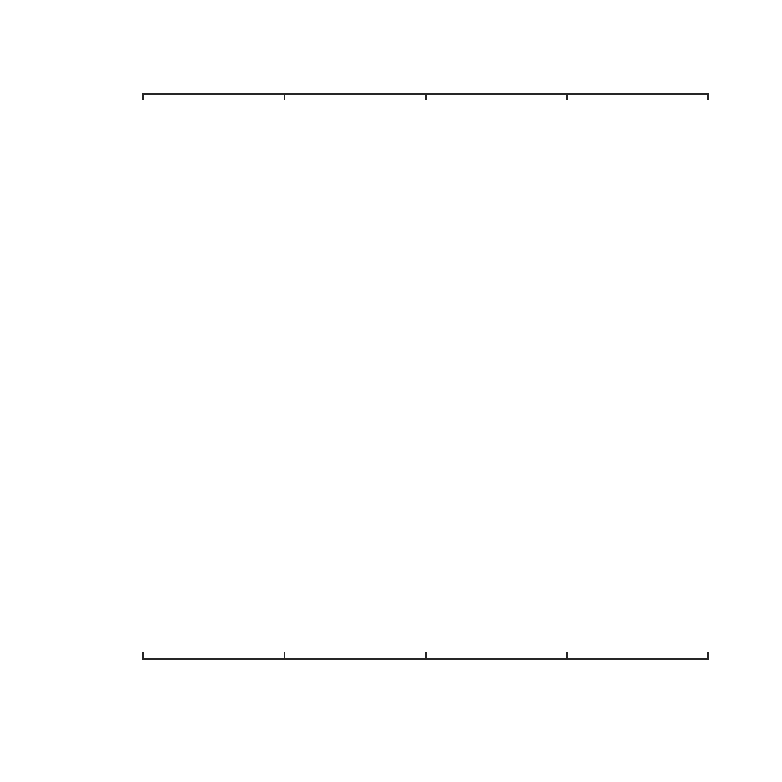}
\caption{}
\end{subfigure}
\caption{Helmholtz projection problem in unit square. (a) Mesh and contours of
  $x$-component of the numerical solution to $\bm{u}$ computed with polynomial
  degree 20. (b) Maximum absolute error in $x-$component $({\color{blue}\circ})$
  and $y-$component $({\color{red}\times})$ of the numerical solution to
  $\bm{u}$ v/s polynomial degree.}
\label{fig:divfreeproj2d}
\end{figure}

Consider $\Omega$ to be a unit square and the triangulation $T_h$ to be a
uniform triangulation with eight elements. The triangulation is shown by the red
lines in figure \ref{fig:divfreeproj2d}a. The function $\bm{g}$ is taken to be
\begin{equation*}
(\cos(2\pi x)\sin(2\pi y)+0.1\cos(2\pi x)\sin(2\pi y),-\sin(2\pi x)\cos(2\pi y)+0.1\sin(2\pi x)\cos(2\pi y)).
\end{equation*}
For this $\bm{g}$, the exact solution $\bm{u}$ and $\lambda$  are
\begin{equation*}
(\cos(2\pi x)\sin(2\pi y),-\sin(2\pi x)\cos(2\pi y)),\text{ and }\frac{0.1}{2\pi}\sin(2\pi x)\sin(2\pi y),
\end{equation*}
respectively. All numerical solutions from degree 0 to 20 are computed. The contours of the $x$-component of the numerical solution to $\bm{u}$ computed for polynomial degree $20$ are shown in figure \ref{fig:divfreeproj2d}a. Figure \ref{fig:divfreeproj2d}b shows the maximum error in the numerical solution to $\bm{u}$ as a function of polynomial degree. This error is the maximum absolute difference between $\bm{u}$ and its numerical solution evaluated at 1600 points in each element. The error decreases exponentially with increasing polynomial degree and reaches $10^{-14}$ for polynomial degree 20.

\begin{figure}[t]
\centering
\begin{subfigure}[b]{0.3\linewidth}
\centering
\adjustbox{max width=\linewidth,trim=0cm 0cm 0.4cm 0.4cm,clip}{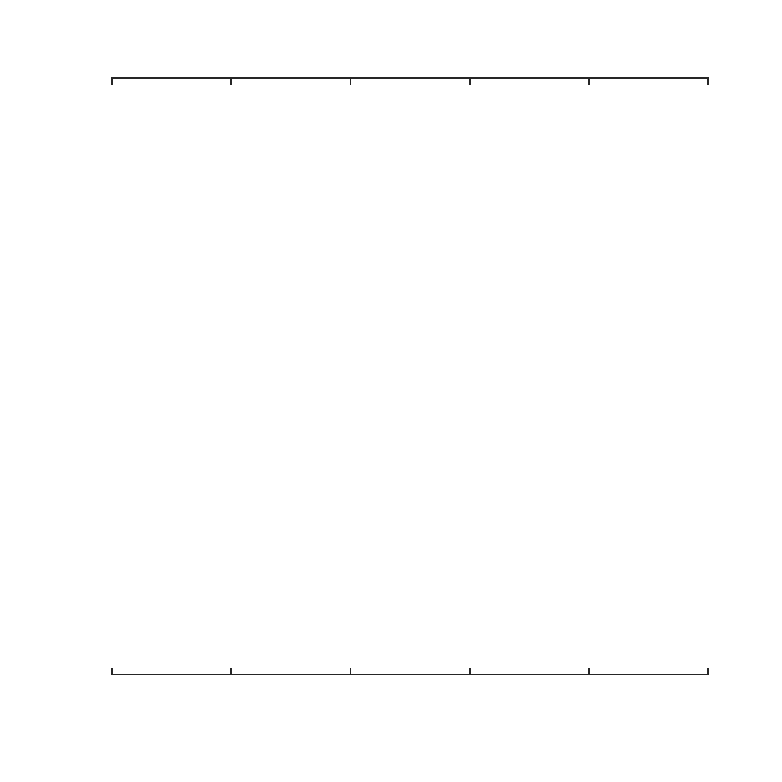}
\caption{}
\end{subfigure}\hfill
\begin{subfigure}[b]{0.35\linewidth}
\centering
\adjustbox{max width=\linewidth,trim=0cm 0.7cm 0.7cm 1cm,clip}{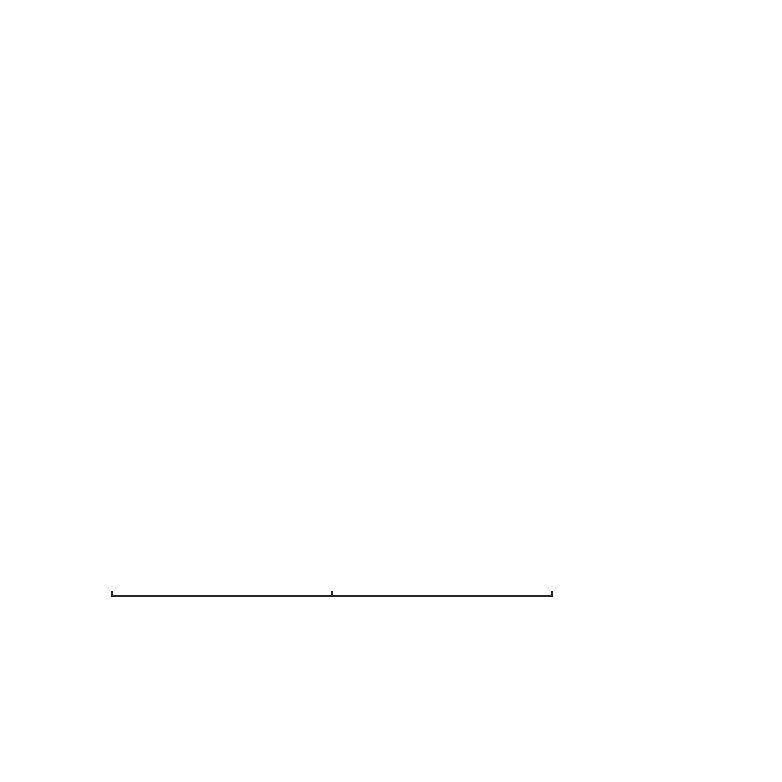}
\caption{}
\end{subfigure}\hfill
\begin{subfigure}[b]{0.3\linewidth}
\centering
\adjustbox{max width=\linewidth,trim=0cm 0cm 0cm 0cm,clip}{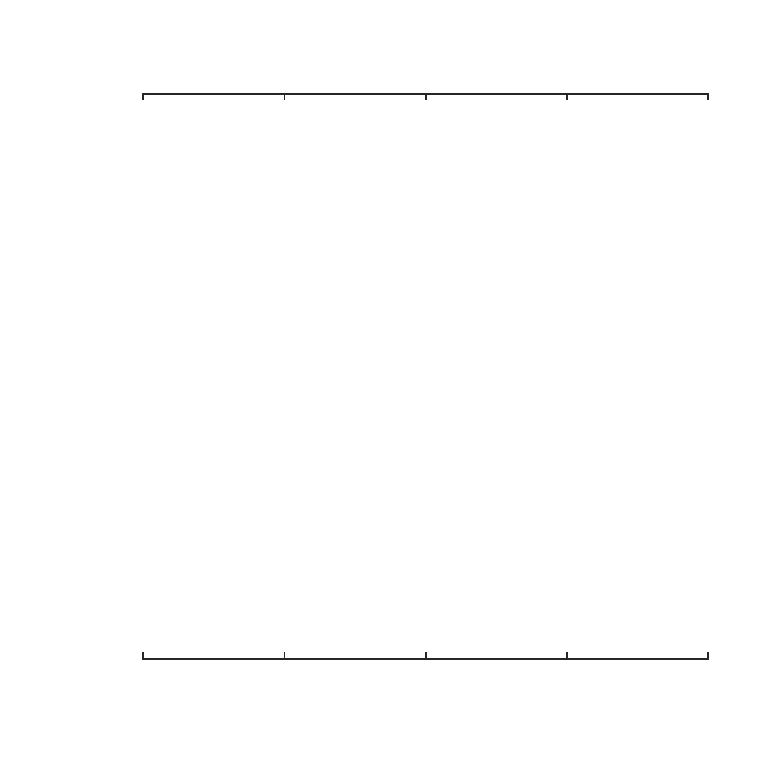}
\caption{}
\end{subfigure}\hfill
\caption{Helmholtz projection problem in a randomly generated two-dimensional domain. (a) Mesh. (b) $x$-component of the numerical solution to $\bm{u}$ computed with polynomial degree 20. (b) Maximum absolute error in $x-$component ($\color{blue}\circ$) and $y-$component ($\color{red}\times$) of the numerical solution to $\bm{u}$ v/s polynomial degree.}
\label{fig:divfreerandproj2d}
\end{figure}

Consider $\Omega$ to be the convex hull of $50$ randomly scattered points in the unit square, and the triangulation $T_h$ to be a Delaunay triangulation of these points. The mesh is shown in figure \ref{fig:divfreerandproj2d}a. The function $\bm{g}$ is taken to be $(\cos(2\pi x)\sin(2\pi y),-\sin(2\pi x)\cos(2\pi y))$. Since $\bm{g}$ is divergence-free, the exact solution $\bm{u}$ equals $\bm{g}$ and $\lambda$ equals zero. All numerical solutions from polynomial degree 0 to 20 are computed. The contours in figure \ref{fig:divfreerandproj2d}b show the $x$-component of the numerical solution to $\bm{u}$ computed with polynomial degree 20. Figure \ref{fig:divfreerandproj2d}c shows the maximum error in the numerical solution to $\bm{u}$ v/s polynomial degree. This error is the maximum absolute difference between $\bm{u}$ and its numerical solution evaluated at 1600 points in each element. The error decreases exponentially and stagnates at machine precision due to round-off error.

\begin{figure}[t]
\centering
\begin{subfigure}[b]{0.32\linewidth}
\centering
\adjustbox{max width=\linewidth,trim=0cm 0cm 0cm 0cm,clip}{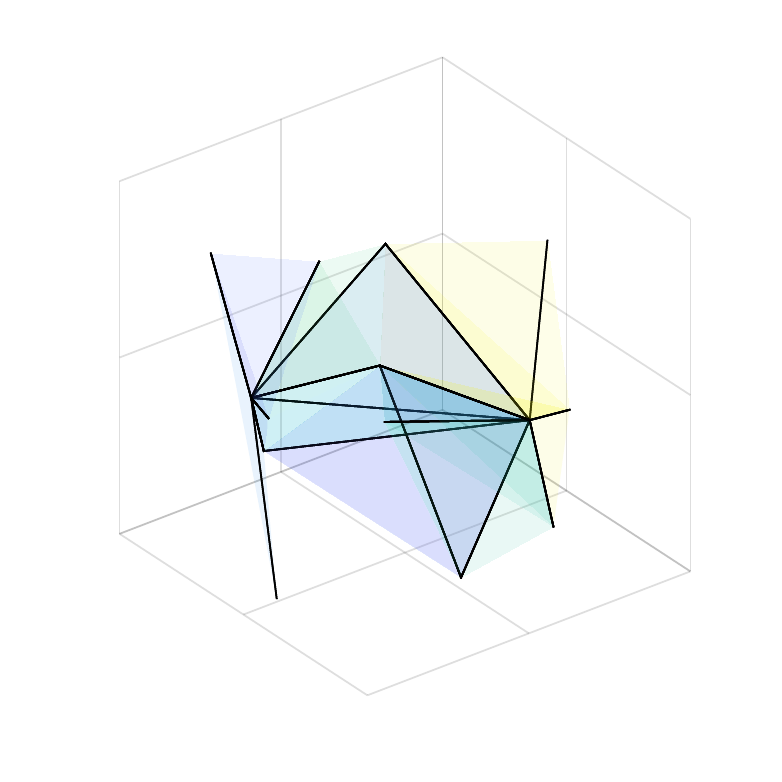}
\caption{}
\end{subfigure}
\begin{subfigure}[b]{0.32\linewidth}
\centering
\adjustbox{max width=\linewidth,trim=0cm 0.7cm 0.7cm 1cm,clip}{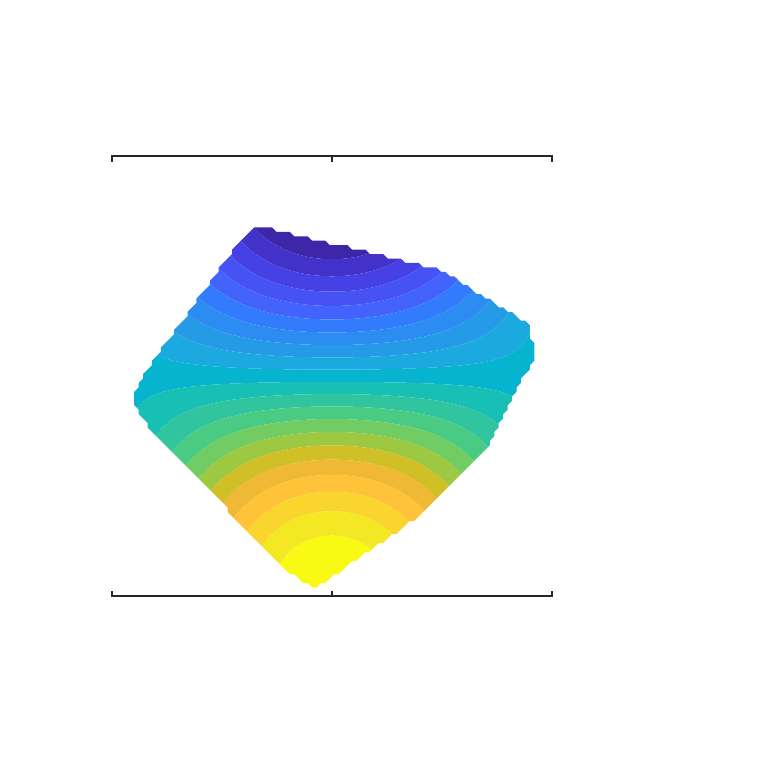}
\caption{}
\end{subfigure}
\begin{subfigure}[b]{0.32\linewidth}
\centering
\adjustbox{max width=\linewidth,trim=0cm 0.7cm 0.7cm 1cm,clip}{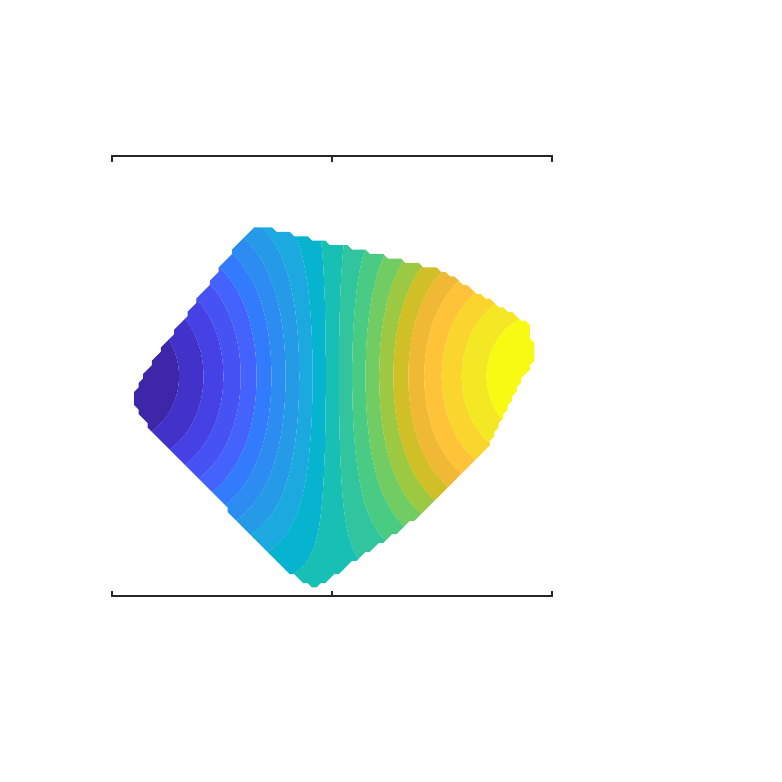}
\caption{}
\end{subfigure}
\begin{subfigure}[b]{0.32\linewidth}
\centering
\adjustbox{max width=\linewidth,trim=0cm 0.7cm 0.6cm 1cm,clip}{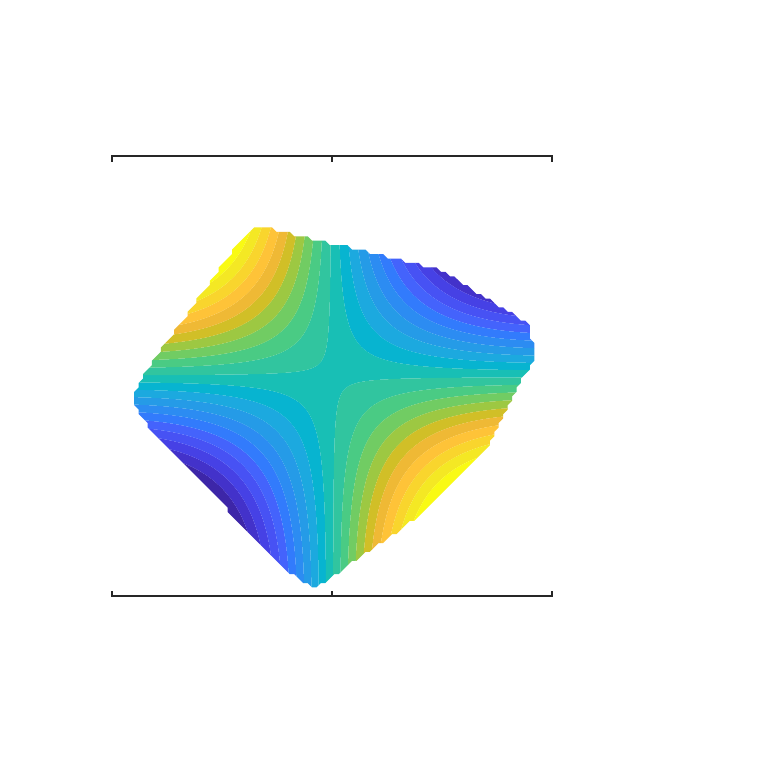}
\caption{}
\end{subfigure}
\begin{subfigure}[b]{0.32\linewidth}
\centering
\adjustbox{max width=\linewidth,trim=0cm 0cm 0cm 0cm,clip}{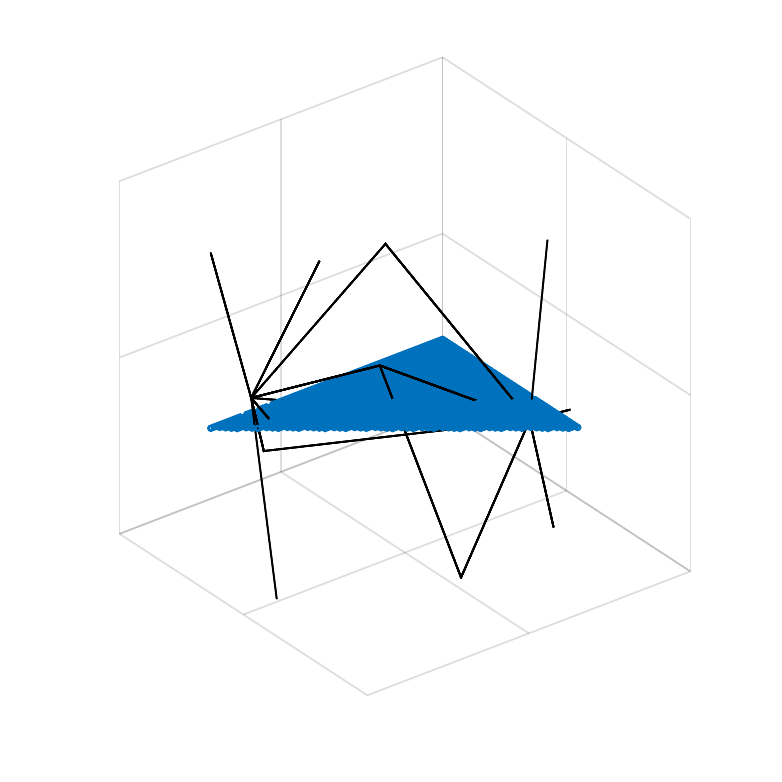}
\caption{}
\end{subfigure}
\begin{subfigure}[b]{0.32\linewidth}
\centering
\adjustbox{max width=\linewidth,trim=0cm 0cm 0cm 0cm,clip}{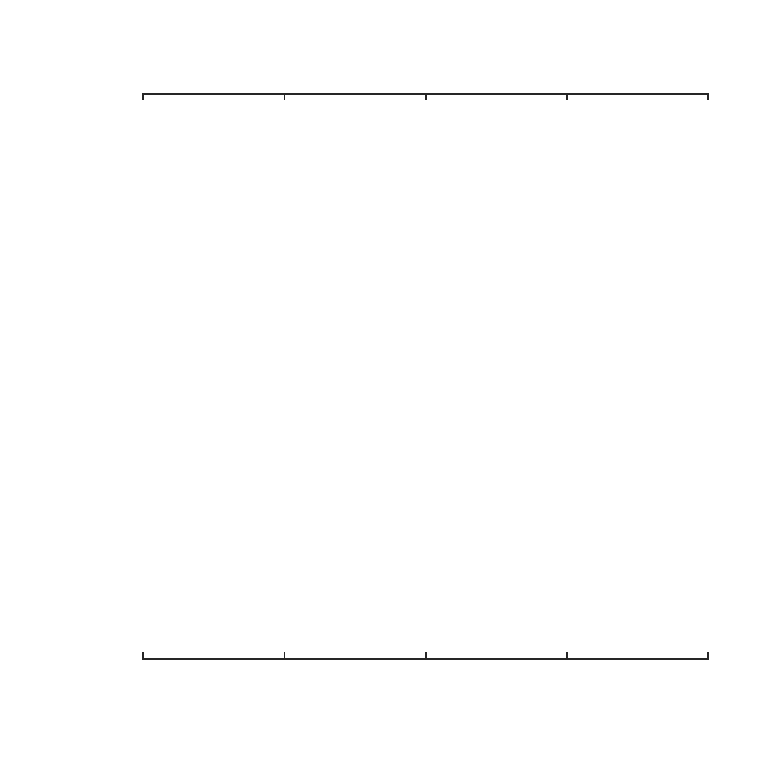}
\caption{}
\end{subfigure}
\caption{Helmholtz projection problem in a randomly generated three dimensional domain. (a) Domain and mesh. (b) Location of z-slices. (c), (d), and (e) show the x-, y-, and z- component of the computed $\vec{u}_h$ on the z-slice. (f) Error in $x-$component ($\color{blue}\circ$), $y-$component ($\color{red}\times$), and $z-$component ($\color{orange}\diamond$) of the numerical solution to $\bm{u}$ v/s polynomial degree.}
\label{fig:divfreeproj3d}
\end{figure}

Consider $\Omega$ to be the three-dimensional convex hull of 20 randomly scattered points in the unit cube and $T_h$ to be a three-dimensional Delaunay triangulation of these points. The domain and the mesh are shown in the figure \ref{fig:divfreeproj3d}a. The function $\bm{g}$ is set to the three dimensional field 
\begin{equation*}
(\sin(\pi x)\cos(\pi y)\cos(\pi z),-0.5\cos(\pi x)\sin(\pi y)\cos(\pi z),-0.5\cos(\pi x)\cos(\pi y)\sin(\pi z)).
\end{equation*}
Since $\bm{g}$ is divergence-free, the exact solution $\bm{u}$ equals $\bm{g}$ and $\lambda$ equals zero. All numerical solutions up to polynomial degree 17 are computed. Figures \ref{fig:divfreeproj3d}b-d show the contours of $x-$, $y-$, and $z-$component of the numerical solution to $\bm{u}$ computed with polynomial degree 10 on a $x-y$ plane located at $z=0.2$ (plane is shown in figure \ref{fig:divfreeproj3d}e). Figure \ref{fig:divfreeproj3d}f shows the maximum error in the numerical solution to $\bm{u}$ v/s polynomial degree. This error is the maximum absolute different between $\bm{u}$ and its numerical solution evaluated at 8000 points in each element. The error decreases exponentially with polynomial degree up to machine epsilon.

\subsection{Laplace problem}
The Laplace problem considered is: Given a domain $\Omega$, and the Dirichlet boundary data $\lambda_D$ on $\partial \Omega$, find $\bm{u}$ and $\lambda$ in $\Omega$ such that
\begin{alignat}{3}
&\bm{u}\quad+&\bm{\nabla}\lambda&=0&&\quad\text{ in }\Omega\\
&\bm{\nabla}\cdot \bm{u}& &=0&&\quad\text{ in }\Omega\nonumber\\
&&\lambda&=\lambda_D&&\quad\text{ on }\partial\Omega\nonumber.
\end{alignat}
Similar to the global divergence-free projection problem, i) we solve the above
equation using our divergence-free basis in place of the usual polynomial basis
to approximate $\bm{u}$ in the hybridized BDM mixed method
\citep{brezzi1985two}, and ii) the left-hand side matrices and the right-hand
side vectors for all polynomial degrees from 0 to the given degree $k$ are
incrementally built.

\begin{figure}[t]
\centering
\begin{subfigure}[b]{0.35\linewidth}
\centering
\adjustbox{max width=\linewidth,trim=0cm 0.7cm 0.7cm 1cm,clip}{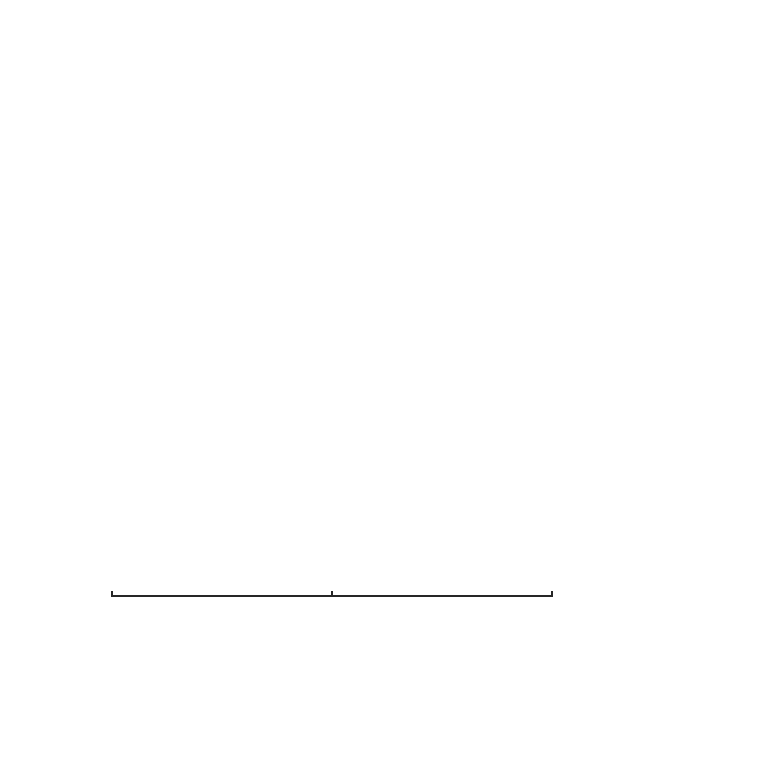}
\caption{}
\end{subfigure}
\begin{subfigure}[b]{0.3\linewidth}
\centering
\adjustbox{max width=\linewidth,trim=0cm 0cm 0cm 0cm,clip}{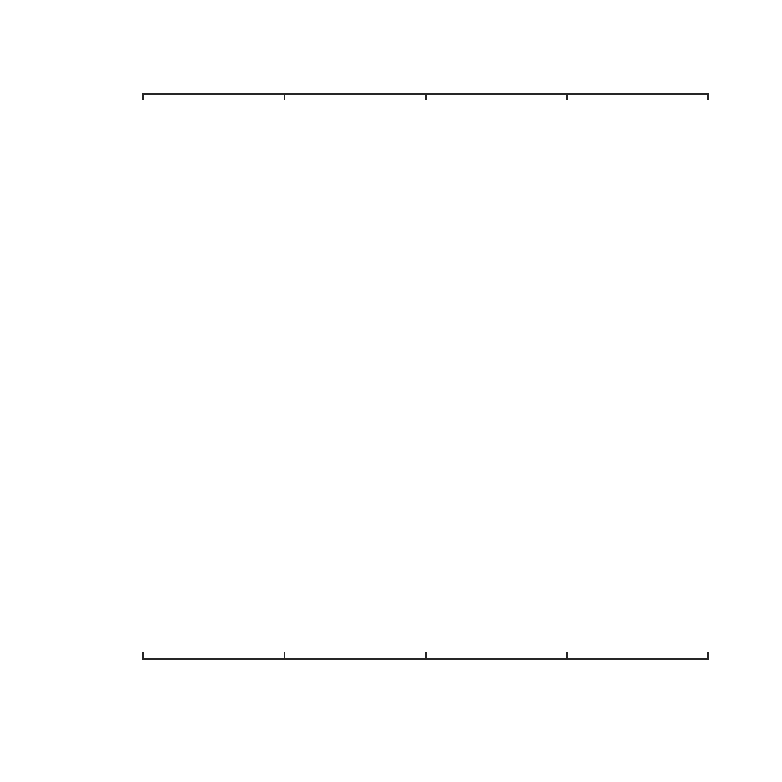}
\caption{}
\end{subfigure}
\caption{Laplace problem. (a) Mesh and $x-$ component of the numerical solution to $\bm{u}$ computed with polynomial degree 20. (b) Error in $x-$component ($\color{blue}\circ$), and $y-$component ($\color{red}\times$) of the numerical solution to $\bm{u}$ v/s polynomial degree.}
\label{fig:lap_2d}
\end{figure}

Consider $\Omega$ to be a unit square and its triangulation $T_h$ to be a
uniform triangulation composed of eight elements. Figure \ref{fig:lap_2d}a shows
the mesh. The boundary data $\lambda_D$ is set using exact solution
$\lambda=\sin(2\pi x)(\cosh(2\pi y)-\sinh(2\pi y))$. Contours of the
$x$-component of the numerical solution to $\bm{u}$ computed with polynomial
degree 20 are shown in figure \ref{fig:lap_2d}a. Figure \ref{fig:lap_2d}b shows
the maximum error in the numerical solution to $\bm{u}$ v/s polynomial
degree. This error is the maximum absolute difference between $\bm{u}$ and its
numerical solution evaluated at 1600 points in each triangle. The error
decreases exponentially with increasing polynomial degree. The error reaches
$10^{-12}$ for degree 15 and then stagnates due to round-off error.

\begin{figure}[t]
\centering
\begin{subfigure}{0.35\linewidth}
\centering
\adjustbox{max width=\linewidth,trim=0cm 0.7cm 0.7cm 1cm,clip}{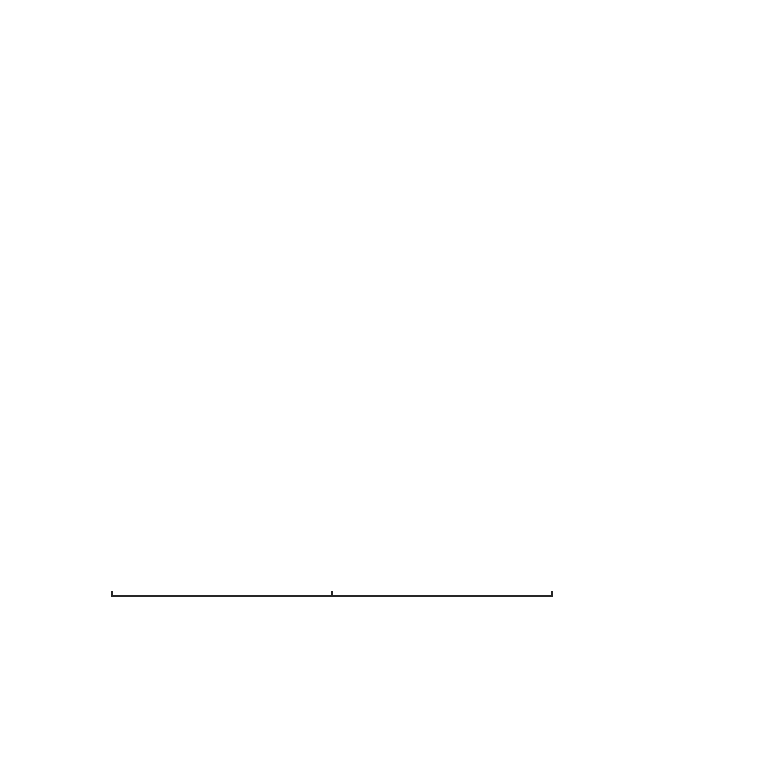}
\caption{}
\end{subfigure}
\begin{subfigure}{0.3\linewidth}
\centering
\adjustbox{max width=\linewidth,trim=0cm 0cm 0cm 0cm,clip}{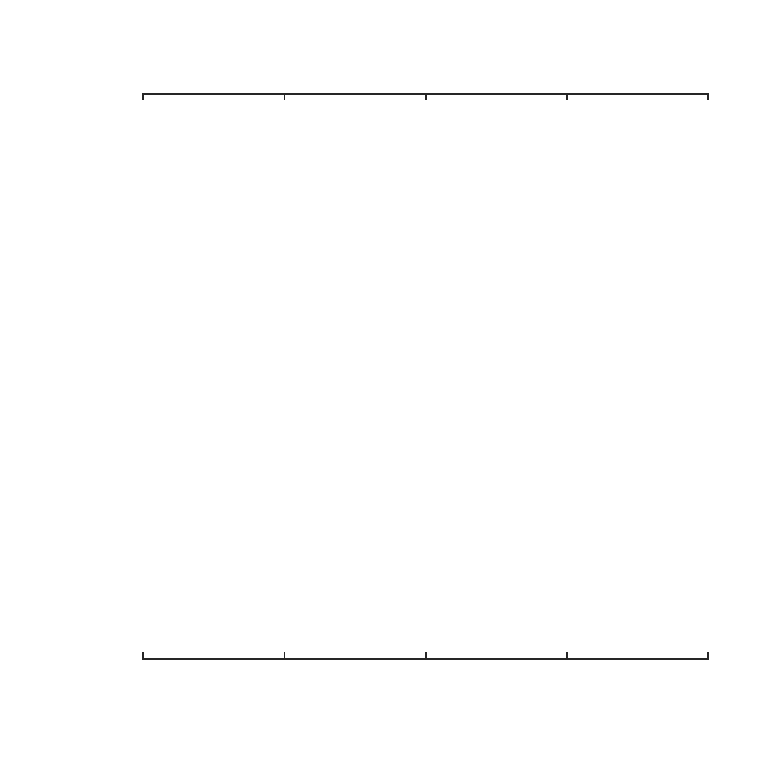}
\caption{}
\end{subfigure}
\caption{Poisson problem. (a) Mesh and $x-$component of the numerical solution to $\bm{u}$ computed with polynomial degree 20. (b) Error in $x-$component ($\color{blue}\circ$), and $y-$component ($\color{red}\times$) of the numerical solution to $\bm{u}$ v/s polynomial degree.}
\label{fig:poiss2d}
\end{figure}

Another Laplace problem considered is the corner singularity problem shown in figure \ref{fig:lap_corner_sing} and its results were discussed in the introduction section. Some important results are reiterated. The numerical solution to $\lambda$ at the corner converges exponentially with increasing polynomial degree and is accurate up to twelve significant digits for polynomial degree eight. It takes just four seconds to compute all numerical solutions from polynomial degree zero to eight. 


\subsection{Poisson problem}
The Poisson problem considered is: Given a domain $\Omega$ and a function $f$ in $\Omega$, find $\bm{u}$ and $\lambda$ such that
\begin{alignat}{3}
&\bm{u}\quad+&\bm{\nabla}\lambda&=0&&\quad\text{ in }\Omega\\
&\bm{\nabla}\cdot \bm{u}& &=f&&\quad\text{ in }\Omega\nonumber\\
&&\lambda&=0&&\quad\text{ on }\partial\Omega\nonumber.
\end{alignat}
The above equations are solved using the hybridized BDM mixed method. In each
element of the mesh, we construct the divergence-free basis functions and use it
to express the portion of $\bm{u}$ that depends on the unknown Lagrange
multiplier on the element faces. For this problem, only the left-hand side
matrix can be incrementally built for all polynomial degrees from zero to the
prescribed degree $k$. The right-hand side vector needs to computed separately
for each degree.

Consider $\Omega$ to be a unit square, and $T_h$ to be a uniform triangulation
composed of eight elements. The mesh is shown in figure \ref{fig:poiss2d}a. The
data $f$ is computed assuming the exact solution $\lambda$ to be
$\sin(2\pi x)\sin(2\pi y)$. Contours of the $x-$component of the numerical
solution to $\bm{u}$ computed with polynomial degree 20 are shown in figure
\ref{fig:poiss2d}a. The error in the numerical solution to $\bm{u}$ is shown in
figure \ref{fig:poiss2d}b as a function of the polynomial degree. This error is
maximum absolute difference between $\bm{u}$ and its numerical solution
evaluated at 1600 points in each triangle. The error decreases exponentially
with increasing polynomial degree. It reaches around $10^{-12}$ for degree 17
and then stagnates due to round-off error.

\section{Summary} \label{sec:summary}

This paper develops a methodology to construct an orthonormal and hierarchical
divergence-free polynomial basis in a simplex (triangles in 2D and tetrahedra in
3D) of arbitrary dimension. At the core of the construction is an Arnoldi-based
procedure that constructs an orthonormal basis for polynomials of degree less
than or equal to $k$ in $d$ dimensions. The generated basis is robust in
finite-precision arithmetic. Using this basis in hybridized mixed methods leads
to fast computation of all numerical solutions from polynomial degree zero to
some given $k$. The orthonormality simplifies the local problem solution. The
hierarchical feature allows the global (and element) matrices and vectors to be
incrementally constructed for all degrees zero to $k$ using the local problem
solution computed just for degree $k$. The constructed basis is applied to solve
Helmholtz, Laplace and Poisson problem in smooth domains and in a domain with
corner singularity. The basis can also be used for efficient numerical solution
of other PDEs such as incompressible Stokes, incompressible Navier-Stokes, and
Maxwell equations.

\section*{Acknowledgements}
This work was supported by the United States Office of Naval Research under
grant N00014-21-1-2454.

\bibliographystyle{cas-model2-names}
\bibliography{papers}

\appendix

\section{The Arnoldi-based procedure generates a basis for
  polynomials} \label{app:proofarpoly}

To show that the Arnoldi-based procedure generates a basis for polynomials, we
can omit the orthogonalizations as they merely combine the polynomials without
adding higher-degree polynomials. The polynomials generated with this omission
are denoted by $m_i$. The algorithm with this omission is: $m_1=1$. For each
degree $j=1,\dots,k$, compute the polynomials $m_{p_{j-1}+1},\dots,m_{p_j}$ of
degree $j$ using the previously computed basis functions $m_1,\dots,m_{p_{j-1}}$
and the algorithm given below:
\begin{algorithmic}[1]
  \State Set $c=p_{j-1}$
  \For{$i=1,\dots,d$}
    \State $j''=C_{d-i}^{j-1+d-i}$
    \For{$j'=1,\dots,j''$}
      \State $m_{c+1}=\xhat_im_{p_{j-1}-j''+j'}$
      \State $c=c+1$
    \EndFor
  \EndFor
\end{algorithmic}

\begin{figure}[h!]
\centering
\begin{subfigure}{0.5\linewidth}
\centering
\begin{tikzpicture}[every node/.style={anchor={center}}]
\matrix(m)[matrix of nodes,column sep=5mm,row sep=5mm]
{$m_1$&$m_2$&$m_4$&$m_7$&$m_{11}$\\
$m_3$&$m_5$&$m_8$&$m_{12}$\\
$m_6$&$m_9$&$m_{13}$\\
$m_{10}$&$m_{14}$\\
$m_{15}$\\};
  { [start chain,every join/.style={->,red}] \chainin (m-1-1);
    \chainin (m-1-2)[join={node[above] {{\small \circled{1}}}}];\chainin (m-1-3)[join={node[above] {{\small \circled{3}}}}];\chainin (m-1-4)[join={node[above] {{\small \circled{6}}}}];\chainin(m-1-5)[join={node[above] {{\small \circled{10}}}}];}
    {[start chain,every join/.style={->,blue}]\chainin(m-1-1);\chainin(m-2-1)[join={node[left] {{\small \circled{2}}}}];\chainin(m-3-1)[join={node[left] {{\small \circled{5}}}}];
    \chainin(m-4-1)[join={node[left] {{\small \circled{9}}}}];
    \chainin(m-5-1)[join={node[left] {{\small \circled{14}}}}];}
   { [start chain,every join/.style={->,red}] \chainin (m-2-1);
    \chainin (m-2-2)[join={node[above] {{\small \circled{4}}}}];\chainin (m-2-3)[join={node[above] {{\small \circled{7}}}}];\chainin (m-2-4)[join={node[above] {{\small \circled{11}}}}];}
    {[start chain,every join/.style={->,red}]\chainin(m-3-1);\chainin(m-3-2)[join={node[above] {{\small \circled{8}}}}];\chainin(m-3-3)[join={node[above] {{\small \circled{12}}}}];}
    {[start chain,every join/.style={->,red}]\chainin(m-4-1);\chainin(m-4-2)[join={node[above] {{\small \circled{13}}}}];}
    {[start chain,every join/.style={->,red}]\chainin(m-5-1);}
    \node also [label=above:($1$)] (m-1-1) {};
    \node also [label=above:($x$)] (m-1-2) {};
    \node also [label=left:($y$)] (m-2-1) {};
    \node also [label=above:($x^2$)] (m-1-3) {}; 
    \node also [label=above:($xy$)] (m-2-2) {};
    \node also [label=left:($y^2$)] (m-3-1) {};
    \node also [label=above:($x^3$)] (m-1-4) {};
    \node also [label=above:($x^2y$)] (m-2-3) {};
    \node also [label=above:($xy^2$)] (m-3-2) {};
    \node also [label=left:($y^3$)] (m-4-1) {};
    \node also [label=above:($x^4$)] (m-1-5) {};
    \node also [label=above:($x^3y$)] (m-2-4) {};
    \node also [label=above:($x^2y^2$)] (m-3-3) {};
    \node also [label=above:($xy^3$)] (m-4-2) {};
    \node also [label=left:($y^4$)] (m-5-1) {};
\end{tikzpicture}
\end{subfigure}
\caption{Illustration of the generation of monomials for $k=4$.}
\label{fig:vis_Pk2d}
\end{figure}
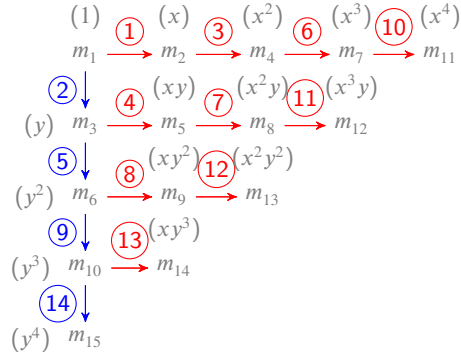

All it remains to prove is that the set of polynomials $m_1,\dots,m_{p_j}$ form
a basis for polynomials of degree less than or equal to $j$. In fact, the $m_i$s
in this set are degree-ordered monomials of degree less than or equal to
$j$. Furthermore, the above algorithm generates the monomials of degree $j$
using the monomials of degree $j-1$ and the coordinate operators. We visually
show this for two dimensions $d=2$. Figure \ref{fig:vis_Pk2d} shows the monomial
generation process for $k=4$ and $d=2$ with the above algorithm. Here, we have
used $x$ and $y$ in place of $x_1$ and $x_2$ for the sake of clarity. The red
and blue arrows denote the action of the $x$-coordinate and $y$-coordinate
operators, respectively. The circled numbers adjacent to each arrow shows the
sequence in which the $m_i$s are generated. Each $m_i$ is equal to the monomial
in parenthesis above or to the left of it. For degree $j$, the monomials are
generated as $x^{j+1-r}y^r=\xhat x^{j-r}y^r$ for $r=1,\dots,j$ and
$y^{j+1}=\yhat y^j$. Similarly, for arbitrary $d$, the monomials of degree $j$
are generated by applying the coordinate operators on the monomials of degree
$j-1$. The first $C_{d-1}^{j-1+d-1}$ monomials of degree $j$ are
generated by applying $\xhat_1$ onto each monomial of degree $j-1$. The next
$C_{d-2}^{j-1+d-2}$ monomials of degree $j$ are generated by applying the
coordinate operator $\xhat_2$ onto the last $C_{d-2}^{j-1+d-2}$ monomials of
degree $j-1$, and so on.

\section{\textbf{Step 1.2} generates an orthonormal basis for divergence-free
  polynomials} \label{app:divfreestep1p2}

The set of divergence-free polynomials is a subset of the set of vector-valued
polynomials. Therefore, each divergence-free basis function can be expanded in
terms of the vector-valued polynomial basis $\{q_r\mbf{e}_i\}$. The
divergence-free polynomial basis function of degree $j$ is
$\bm{\varphi}_{\ell}=\sum_{i=1}^d\sum_{r=1}^{p_j}N_{(i-1)p+r,\ell}q_r\bm{e}_i$. For
degree $j$ and $j-1$, the number of divergence-free basis functions are $n_j$
and $n_{j-1}$, respectively. Therefore, the loop for degree $j$ should add
$n_j-n_{j-1}$ functions, i.e.,
$\bm{\varphi}_{n_{j-1}+1},\dots,\bm{\varphi}_{n_j}$. These functions must be
divergence-free. Hence, the divergence-free condition (i). The divergence of a
degree $j$ vector-valued polynomial is another polynomial of degree $j-1$. The
condition (i) is an integral condition that requires this degree $j-1$
polynomial to be zero. It does so by requiring its projection to be zero along
each orthonormal polynomial $q_i$ of degree less than or equal to $j-1$. We use
an integral condition for numerical stability reasons in the presence of
finite-precision arithmetic. Using a pointwise imposition of the divergence-free
constraint lead to large amplification of machine precision error at high
polynomial degrees.

The dimension of the vector-valued polynomial basis $\{q_r\mbf{e}_i\}$ for
$r=1,\dots,p_j$ and $i=1,\dots,d$ is $dC_d^{j+d}$. The divergence-free condition
imposes $C_d^{j-1+d}$ constraints. This gives us a total dimension of
$n_j=dC_d^{j+d}-C_d^{j-1+d}$. To determine the new $n_j-n_{j-1}$ functions, we
need certain other constraints. We require the new functions to have no
component along the divergence-free basis funcitons of degree less than or equal
to $j-1$, i.e.,
$\int_{\Omegahat}\bm{\varphi}_{\ell}\cdot\bm{\varphi}_{\ell'}\,d\Omegahat=0$ for
$\ell=n_{j-1}+1,\dots,n_j$ and $\ell'=1,\dots,n_{j-1}$. Combining this with the
requirement that the new $n_j-n_{j-1}$ must be orthonormal amongst each other
yields condition (ii). Lastly, we would not want the new basis functions to be
linearly dependent which would reduce the total dimension from being equal to
$n_j$. Hence, we have condition (iii).

\section{Quadrature rule} \label{app:quadrule}


\begin{figure}
\centering
\begin{subfigure}{0.32\linewidth}
\centering
\adjustbox{max width=\linewidth,trim=0cm 0cm 0cm 0cm,clip}{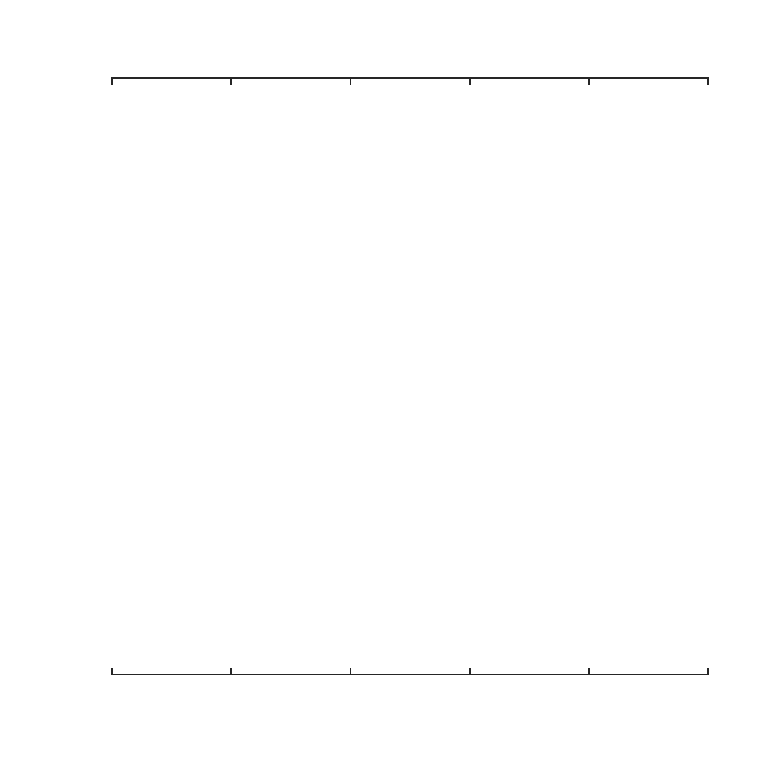}
\caption{}
\end{subfigure}
\begin{subfigure}{0.32\linewidth}
\centering
\adjustbox{max width=\linewidth,trim=0cm 0cm 0cm 0cm,clip}{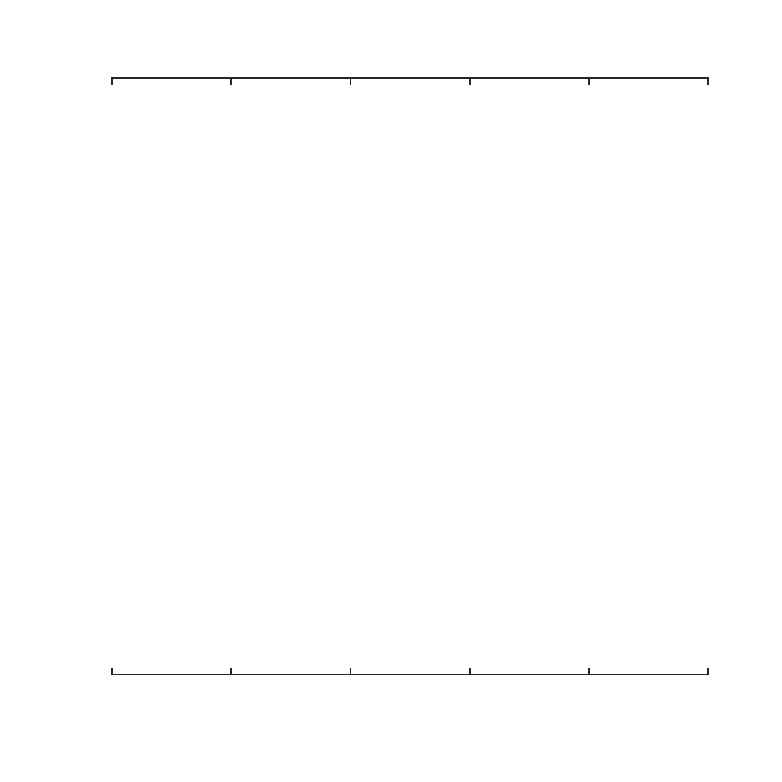}
\caption{}
\end{subfigure}
\caption{(a) Quadrature points to exactly integrate a polynomial of degree 4 in the unit triangle. (b) Corresponding coordinates of the quadrature points in the unit square.}
\label{fig:quadpts}
\end{figure}

We first derive the quadrature rule for a unit triangle. Consider the following
integral over a unit triangle: $g=\int_0^1\int_0^{1-x}f(x,y)\,dy\,dx$. Map the
unit triangle to a unit square. Define $x=\zeta$ and $y=(1-\zeta)\eta$. Here,
$(\zeta,\eta)$ is the coordinate of a point in the unit square, and $(x,y)$ is
the image of this point in the unit triangle. The determinant of the Jacobian of
this mapping is $(1-\zeta)$. Therefore, the integral transforms to
$g=\int_0^1\int_0^1f(x(\zeta,\eta),y(\zeta,\eta))(1-\zeta)\,d\zeta\,d\eta$ in
the unit square. Consider $f$ to be a polynomial of degree $2k$ in $x$ and $y$,
then $f(x(\zeta,\eta),y(\zeta,\eta))$ is a polynomial in $\zeta$ and $\eta$ with
each exponent less than or equal to $2k$. To integrate exactly along $\eta$, use
$(k+1)$ Gauss-Legendre quadrature points along $\eta$. Along $\zeta$, note that
there is the weight $(1-\zeta)$. Therefore, to integral exactly along $\zeta$,
use $(k+1)$ Gauss-Jacobi points that correspond to the weight
$(1-\zeta)$. Therefore, the quadrature rule to exactly integrate a polynomial of
degree $2k$ is
$\sum_{i=1}^{k+1}\sum_{j=1}^{k+1}f(x_{i,j},y_{i,j})w_{i,j}$. Here,
$x_{i,j}=\zeta_i$, $y_{i,j}=(1-\zeta_i)\eta_j$, and
$w_{i,j}=w_i^{\zeta} w_j^{\eta}$. $\{\zeta_i\}_{i=1}^{k+1}$ and
$\{w_i^{\zeta}\}_{i=1}^{k+1}$ are the Gauss-Jacobi quadrature points and weights
with the weight function $(1-\zeta)$, respectively. $\{\eta_i\}_{i=1}^{k+1}$ and
$\{w_i^{\eta}\}_{i=1}^{k+1}$ are the Gauss-Legendre quadrature points and
weights, respectively. Figure \ref{fig:quadpts}a shows the quadrature points to
exactly integrate a polynomial of degree 4 in the unit triangle. The
corresponding points in the unit square are shown in figure
\ref{fig:quadpts}b. Similarly, in arbitrary dimension $d$, the general integral
over the unit simplex is:
\begin{equation*}
  g=\int_{0}^1\int_0^{1-x_1}\int_0^{1-x_1-x_2}\dots\int_0^{1-x_1-\dots-x_{d-1}}f\,dx_d\,\dots dx_1.
\end{equation*}
The mapping from the unit hypercube to the unit simplex $\bm{x}(\bm{\zeta})$ is:
\begin{equation*}
  x_1=\zeta_1,\,x_2=(1-\zeta_1)\zeta_2,\dots,\,x_d=(1-\zeta_1)\dots(1-\zeta_{d-1})\zeta_d.
\end{equation*}
The transformed integral is:
\begin{equation*}
  g=\int_{0}^1\dots\int_0^{1}f(\bm{x}(\bm{\zeta}))(1-\zeta_1)^{(d-1)}\dots(1-\zeta_{d-1})\,d\zeta_d\,\dots d\zeta_1.
\end{equation*}
The quadrature rule is:
\begin{equation*}
  g=\sum_{i_1=1}^{k+1}\dots\sum_{i_d=1}^{k+1}f(x_{1_{i_1,\dots,i_d}},\dots,x_{d_{i_1,\dots,i_d}})w_{i_1,\dots,i_d}.
\end{equation*}
Here, the quadrature point
$\bm{x}_{i_1,\dots,i_d}=\bm{x}(\zeta_{1_{i_1}},\dots,\zeta_{d_{i_d}})$ and
$\{\zeta_{j_i}\}_{i=1}^{k+1}$ are the Gauss-Jacobi quadrature points with weight
function $(1-\zeta_j)^{d-j}$. The quadrature weight
$w_{i_1,\dots,i_d}=w_{1_{i_1}}\dots w_{d_{i_d}}$, where
$\{w_{j_i}\}_{i=1}^{k+1}$ are the corresponding Gauss-Jacobi quadrature weights
along the $j^{th}$ dimension. The quadrature points and weights are rearranged
into a matrix $x$ and a vector $w$, respectively. $x$ is a matrix of size
$(k+1)^d\times d$. The $i^{th}$ component of the coordinate vector of each point
is stored in the $i^{th}$ column of $x$. $w$ is a vector of size $(k+1)^d$.

\section{Construction of derivative matrices} \label{app:dermat}

The derivative matrices yield the partial derivative of the function at the
quadrature points using the values of the function at the same points. We first
derive these matrices for two dimensions. Observe that in the unit square
(figure \ref{fig:quadpts}b), the quadrature points form a Cartesian grid. The
derivative matrices for this Cartesian grid can be computed using the
one-dimensional derivative matrices and taking its kronecker tensor product with
the identity matrix. The derivative matrix in the unit triangle can then be
obtained using the chain rule:
$\partial f/\partial x=\partial f/\partial \zeta\partial \zeta /\partial
x+\partial f/\partial \eta\partial \eta/\partial x$ and
$\partial f/\partial y=\partial f/\partial \zeta\partial \zeta /\partial
y+\partial f/\partial \eta\partial \eta/\partial y$. Denote the derivative
matrix along the $\zeta$ direction by $DZ$. It is the kronecker tensor product
$DZ=kron(I,DZ1)$. Here, $I$ is the $(k+1)\times (k+1)$ identity matrix and $DZ1$
is the one-dimensional differentiation matrix along $\zeta$. To construct $DZ1$,
we use the barycentric Lagrange-based procedure of
\cite{berrut2004barycentric}. $DZ1$ is a $(k+1)\times (k+1)$ matrix given by
$DZ1_{i,j}=(\lambda_j/\lambda_i)1/(\zeta_i-\zeta_j)$ for $i\neq j$ and
$DZ1_{i,j}=-\sum_{i\neq j}DZ1_{i,j}$. Here, $\zeta_i$ is the $i^{th}$ quadrature
point along the $\zeta$ direction. $\lambda_i$ is the $i^{th}$ barycentric
weight given by $\lambda_i=1/\prod_{j\neq i}(\zeta_i-\zeta_j)$. Similarly, the
derivative matrix along $\eta$ direction can be constructed. It is
$DN=kron(DN1,I)$, where $DN1$ is the one-dimensional differentiation matrix
along $\eta$ given by $DN1_{i,j}=(\gamma_j/\gamma_i)1/(\eta_i-\eta_j)$ for $i\neq j$
and $DN1_{i,j}=-\sum_{i\neq j}DN1_{i,j}$. Here, $\eta_i$ is the $i^{th}$
quadrature point along the $\eta$ direction. $\gamma_i$ is the $i^{th}$
barycentric weight given by $\gamma_i=1/\prod_{j\neq i}(\eta_i-\eta_j)$. In the
unit simplex, the derivative matrices are then $DX=DZ+(y/(1-x)^2).^*DN$ and
$DY=(1/(1-x)).^*DN$, where the $.^*$ operator follows the MATLAB notation.

This procedure to construct the derivative matrices can be extended to an
arbitrary dimension $d$. Denote the matrix of quadrature points by $x$. Its size
is $(k+1)^d\times d$. Let $\{\zeta^{(i)}_j\}_{j=1}^{k+1}$ denote the
Gauss-Jacobi quadrature points in the unit hypercube along the $i^{th}$
direction. The algorithm is:
\begin{algorithmic}[1]
  \For{$i=1,\dots,d$}\Comment{Barycentric weights}
    \For{$j=1,\dots,k+1$}
      \State $\lambda^{(i)}_j=1/\prod_{j\neq \ell}(\zeta^{(i)}_j-\zeta^{(i)}_{\ell})$
    \EndFor
  \EndFor
  \For{$i=1,\dots,d$}\Comment{One-dimensional derivative matrices in the unit hypercube}
    \For{$m=1,\dots,k+1$}
      \For{$n=1,\dots,k+1$}
        \If{$m\neq n$}
          \State $DZ1^{(i)}_{m,n}=\lambda^{(i)}_n/\lambda^{(i)}_m1/(\zeta^{(i)}_m-\zeta^{(i)}_n)$
        \Else
          \State $DZ1^{(i)}_{m,m}=-\sum_{j=1,j\neq m}^{k+1}\lambda^{(i)}_j/\lambda^{(i)}_m1/(\zeta^{(i)}_m-\zeta^{(i)}_j)$
        \EndIf
      \EndFor
    \EndFor
  \EndFor
  \State $DZ^{(1)}=DZ1^{(1)}$ \Comment{Multi-dimensional derivative matrices in the unit hypercube}
  \For{$i=2,\dots,d$}
    \State $DZ^{(i)}=eye(k+1)$
  \EndFor
  \For{$i=1,\dots,d$}
    \For{$j=2,\dots,d$}
      \If{$j==i$}
        \State $DZ^{(i)}=kron(DZ1^{(i)},DZ^{(i)})$
      \Else
        \State $DZ^{(i)}=kron(eye(k+1),DZ^{(i)})$
      \EndIf
    \EndFor
  \EndFor
  \State $DX^{(1)}=DZ^{(1)}$ \Comment{Transforming the derivative matrices from unit hypercube to unit simplex}
  \For{$i=2,\dots,d$}
    \State $DX^{(i)}=zeros((k+1)^d)$
  \EndFor
  \State $tt1=ones((k+1)^d)$
  \For{$j=2,\dots,d$}
    \State $tt1=tt1-x(:,j-1)$; $tt2=1./tt1$; $tt3=x(:,j)./tt1^2$
    \For{$i=1,\dots,j-1$}
      \State $DX^{(i)}=DX^{(i)}+DZ^{(j)}.^*tt3$
    \EndFor
    \State $DX^{(j)}=DX^{(j)}+DR^{(j)}.^*tt2$
  \EndFor
\end{algorithmic}

\section{Deriving the algorithm in section \ref{subsubsec:1p2} from the
algorithm in \textbf{step 1.2} in section
\ref{sec:divfreebasesconst}} \label{app:1p2}

Consider the three conditions in \textbf{step 1.2} of section
\ref{sec:divfreebasesconst}. Substituting the expression for
$\bm{\varphi}_{\ell}$ into condition (i) yields
\begin{equation*}
  \sum_{i=1}^d\sum_{r=1}^{p_j}\left(\int_{\Omegahat}q_s\frac{\partial q_r}{\partial x_i}\,d\Omegahat\right)N_{(i-1)p+r,\ell}=0\text{ for }\ell=n_{j-1}+1,\dots,n_j\text{ and }s=1,\dots,p_{j-1}.
\end{equation*}
These are linear constraints imposed on $\{N_{(i-1)p+r,\ell}\}$. Define the
coefficients of the constraint to be
$C_{s,(i-1)p+r}=\left(\int_{\Omegahat}q_s\frac{\partial q_r}{\partial
    x_i}\,d\Omegahat\right)$. The coefficients $C_{s,(i-1)p+r}$ are stored in
the divergence-free constraint matrix $C$. Using the orthonormality of the
$q_i$s, condition (ii) can be shown to be equivalent to
$\sum_{i=1}^d\sum_{r=1}^{p_j}N_{(i-1)p+r,\ell'}N_{(i-1)p+r,\ell}=\delta_{\ell
  \ell'}$ for $\ell=n_{j-1}+1,\dots,n_j,\ell'=1,\dots,n_j$. Since the $q_i$s are
linearly independent, condition (iii) is equivalent to requiring the rank of the
submatrix of $N$ formed by the $n_{j-1}+1$ to $n_j$ columns to be
$n_j-n_{j-1}$. Computing the coefficients $\{N_{(i-1)p+r,\ell}\}$ that satisfy
the above conditions for each polynomial degree $j$ is equivalent to finding an
orthonormal basis for the null-space of the augmented matrix:
\begin{equation*}
\begin{bmatrix}C_{1:p_{j-1},ii}\\ N^T_{ii,1:n_{j-1}}\end{bmatrix}.
\end{equation*}
Here, $ii$ is the index vector of column indices of $C$ (or row indices of $N$)
that correspond to degree less than or equal to $j$ constructed as:
\begin{algorithmic}[1]
  \State $ii=[]$
  \For{$i=1,\dots,d$}
    \State $ii=[ii,\,(i-1)p+1:(i-1)p+p_j]$
  \EndFor
\end{algorithmic}
Replacing symbolic polynomials with vectors storing its value at the quadrature
points and continous inner-products with equivalent discrete inner-products
yields the algorithm in \ref{subsubsec:1p2}.

\section{Deriving the algorithm in section \ref{subsubsec:2p1} from the
algorithm in \textbf{step 2.1} in section
\ref{sec:divfreebasesconst}} \label{app:2p1}

Using the orthonormality of the $q_i$s, condition (ii) simplifies to
$\sum_{i=1}^d\sum_{r=1}^{p_j}N_{(i-1)p+r,\ell}^{(e)}N_{(i-1)p+r,\ell'}^{(e)}=\delta_{\ell\ell'}$
for $\ell$ and $\ell'=n_{j-1}+1,\dots,n_j$. Since the $q_i$s are linearly
independent, condition (iii) simplifies to enforcing the rank of the submatrix
of $N^{(e)}$ formed by the columns from $n_{j-1}+1$ to $n_j$ to be
$n_{j-1}-n_j$. To compute the coefficients that satisfy the divergence-free
requirement in condition (i), we use the divergence-free basis functions
$\bm{\varphi}_{n_{j-1}+1},\dots,\bm{\varphi}_{n_j}$ that were computed in the
reference element. The expression for these functions in terms of the
orthonormal polynomials is
$\bm{\varphi}_{\ell}=\sum_{i=1}^d\sum_{r=1}^{p_j}N_{(i-1)p+r,\ell}q_r\bm{e}_i,\text{
  for }\ell=n_{j-1}+1,\dots,n_j$. The coefficients $N_{(i-1)p+r,\ell}$ were
computed such that they satisfy the divergence-free condition in the reference
element:
$\sum_{i=1}^d\sum_{r=1}^{p_j}{\partial q_r}/{\partial
  x_i}N_{(i-1)p+r,\ell}=0$. The gradient ${\partial q_r}/{\partial x_i}$ in the
reference element can be expanded in terms of the gradient in the current
element as
${\partial q_r}/{\partial x_i}=\sum_{m=1}^d{\partial x_m^{(e)}}/{\partial
  x_i}{\partial q_r}/{\partial x_m^{(e)}}$. Substituting this into the
divergence-free condition gives
$\sum_{i=1}^d\sum_{r=1}^{p_j}\sum_{m=1}^d{\partial x_m^{(e)}}/{\partial
  x_i}{\partial q_r}/{\partial x_m^{(e)}}N_{(i-1)p+r,\ell}=0$. Swapping the
order of summation over $i$ and $m$, and rearranging yields
$\sum_{m=1}^d\sum_{r=1}^{p_j}\left[\sum_{i=1}^d{\partial x_m^{(e)}}/{\partial
    x_i}N_{(i-1)p+r,\ell}\right]{\partial q_r}/{\partial
  x_m^{(e)}}=0$. Interchanging the index $i$ to $m$ and $m$ to $i$ gives
$\sum_{i=1}^d\sum_{r=1}^{p_j}\left[\sum_{m=1}^d{\partial x_i^{(e)}}/{\partial
    x_m}N_{(m-1)p+r,\ell}\right]{\partial q_r}/{\partial
  x_i^{(e)}}=0$. Comparing the above equation to the divergence-free condition
(i), we deduce that setting $N_{(i-1)p+r,\ell}^{(e)}$ to a linear combination of
the form
$\sum_{\ell'=n_{j-1}+1}^{n_j}A_{\ell',\ell}\left[\sum_{m=1}^d{\partial
    x_i^{(e)}}/{\partial x_m}N_{(m-1)p+r,\ell'}\right]$ will satisfy the
divergence-free condition, for any set of coefficients $A_{\ell',\ell}$ that
form a full rank matrix $A_{n_{j-1}+1:n_j, n_{j-1}+1:n_j}$. The coefficients
$A_{\ell',\ell}$ are implicitly chosen as follows such that the resulting
$N_{(i-1)p+r,\ell}^{(e)}$s satisfy the orthonormality condition (ii). The linear
combinations
$\bar{N}_{(i-1)p_j+r,\ell}=\sum_{m=1}^d\partial x_i^{(e)}/\partial
x_mN_{(m-1)p+r,\ell}$ are computed for $i=1,\dots,d$, $r=1,\dots,p_j$, and
$\ell=n_{j-1}+1,\dots,n_j$. These linear combinations are orthogonalized against
the previously computed coefficients $N_{(i-1)p+r,\ell}^{(e)}$ such that new
$\bar{N}_{(i-1)p_j+r,\ell}$s satisfy
$\sum_{i=1}^d\sum_{r=1}^{p_j}\bar{N}_{(i-1)p_j+r,\ell}N_{(i-1)p+r,\ell'}^{(e)}=0$
for $\ell'=1,\dots,n_{j-1}$. These $\bar{N}_{(i-1)p_j+r,\ell}$s are finally
orthonormalized amongst each other to yield the required coefficients
${N}_{(i-1)p+r,\ell}^{(e)}$s that satisfy
$\sum_{i=1}^d\sum_{r=1}^{p_j}{N}_{(i-1)p+r,\ell}^{(e)}N_{(i-1)p+r,\ell'}^{(e)}=0$
for $\ell$ and $\ell'$ equals $n_{j-1}+1,\dots,n_j$. The computed coefficients
${N}_{(i-1)p+r,\ell}^{(e)}$s also satisfy the rank requirement condition (iii)
because the matrix formed by the linear combinations
$\bar{N}_{(i-1)p_j+r,\ell}$s has rank $n_{j}-n_{j-1}$ and the subsequent
orthogonalizations and orthonormalizations that yield the $n_{j-1}$ to $n_j$
columns of $N^{(e)}$ do not modify this rank. This yields the algorithm in
section \ref{subsubsec:2p1}.

\section{Monomial divergence-free basis functions} \label{app:mondivfreebf}

Figure \ref{fig:mondivfreebf} shows the MATLAB code to construct the monomial divergence-free basis for arbitrary polynomial degree and spatial dimension.

\begin{figure}
\centering
\begin{minipage}{\columnwidth}
\begin{lstlisting}
function [bf] = mondivfreebf(k,d)
%%%% function [bf] = mondivfreebf(k,d)
%%%% Inputs: k - polynomial degree, d - spatial dimension.
%%%% Outputs: bf - function handle to the divergence-free basis funcitons.  
  sz=d*nchoosek(k+d,d)-nchoosek(k-1+d,d); alpha=zeros(d,sz); l=zeros(d,d,sz); t=0; bf=cell(sz,1);
%%%%--------------------Loop over polynomial degree.
  for ktil=0:k
%%%%--------------------Construct monomial basis of homogeneous polynomials of degree ktil.    
    [talpha,tl]=divfree_subspace_of_ptilk(ktil,d); tsz=size(talpha,2);
%%%%--------------------Get function handles.
    for i=1:tsz, tta=talpha(:,i); ttl=tl(:,:,i); bf{t+i}=@(x)divfree_bf_gen_func(x,tta,ttl); end
    t=t+tsz;
  end
function [alpha,l] = divfree_subspace_of_ptilk (k,d)
  if(k==0), alpha=eye(d); l=zeros(d,d,d); return; end
  szptilk=nchoosek(k+d-1,d-1); szptilkd=d*szptilk; szdivc=nchoosek(k+d-2,d-1); szdivptilkd=szptilkd-szdivc;
  alpha=zeros(d,szdivptilkd); l=zeros(d,d,szdivptilkd);
  mnptilk=monomials_of_ptilk(k,d); mnptilkm1=monomials_of_ptilk(k-1,d); eqnmapkm1=mnexp_to_eqnmap(mnptilkm1,k-1,d);
  coeff=zeros(d,szdivc); idcoeff=zeros(d,szdivc); t=0;
  for i=1:d
    for n=1:szptilk
      % Check for zero divergence.      
      tt1=mnptilk(i,n);
      if(tt1==0), t=t+1; alpha(i,t)=1d0; l(:,i,t)=mnptilk(:,n);
      else % We have non-zero divergence.
        talpha=tt1; tl=mnptilk(:,n); tl(i)=tl(i)-1;
        tl1d=mnexp_to_1d(tl,k-1,d); eqnno=eqnmapkm1(tl1d); coeff(i,eqnno)=talpha; idcoeff(i,eqnno)=n;
      end
    end
  end
  % Construct the remaining basis functions.  
  tt2=zeros(d-1,1);
  for n=1:szdivc
    % Solve for the first coefficient in terms of the remaining coefficients.
    for i=1:d-1, tt2(i)=-coeff(i+1,n)/coeff(1,n); end
    % Construct the remaining divergence-free basis functions using these monomials.
    for i=2:d
      t=t+1; alpha(1,t)=tt2(i-1); l(:,1,t)=mnptilk(:,idcoeff(1,n)); alpha(i,t)=1d0;
      l(:,i,t)=mnptilk(:,idcoeff(i,n));
    end
  end  
function [v] = monomials_of_ptilk (k,d)
  tmp=nchoosek(1:(k+d-1),d-1); tmp=tmp'; szptilk=size(tmp,2); v=zeros(d,szptilk);
  for n=1:szptilk
    v(1,n)=tmp(1,n)-1;
    for i=2:d-1
      v(i,n)=tmp(i,n)-tmp(i-1,n)-1;
    end
    v(d,n)=k+d-1-tmp(d-1,n);
  end
function [v] = mnexp_to_eqnmap (mn,k,d)
  v=zeros((k+1)^d,1);
  sz=size(mn,2);
  for n=1:sz, tmp=mnexp_to_1d(mn(:,n),k,d); v(tmp)=n; end
function [v] = mnexp_to_1d (mn,k,d)
  v=mn(1)+1;
  for i=2:d, v=v+(mn(i))*((k+1)^(i-1)); end  
function [bf] = divfree_bf_gen_func(x,talpha,tl)
  [m n]=size(x); bf=zeros([m n]);
  for i=1:n, bf(:,i)=talpha.*(prod(x(:,i).^tl(:,:))'); end
\end{lstlisting}
\end{minipage}
\caption{\texttt{mondivfreebf.m} - MATLAB code to construct the monomial divergence-free basis.}
\label{fig:mondivfreebf}
\end{figure}

\section{MATLAB implementation of step 1} \label{app:matlabstep1}

The MATLAB function \texttt{ardivfreebfref} in figure \ref{fig:ardivfreebfref}
shows our MATLAB implementation of \textbf{step 1}. The required quadrature
points and weights are constructed in lines 9-12. Lines 31-38 show the
operations of the Arnoldi-based process that correspond to the $j^{th}$
polynomial degree. The function \texttt{GaussJacobi} in line 9 is an external
function that yields the one-dimensional quadrature points \texttt{x1} and
weights \texttt{w1} for the $(k+1)$-point Gauss-Jacobi quadrature with weight
function $(1-x)^{\texttt{d-i}}$ in the interval $[-1,+1]$. We use the
\texttt{GuassJacobi} function from
\url{https://www.math.umd.edu/~petersd/460/GaussJacobi.m}.  Lines 14-23 show the
construction of the derivative matrices \texttt{DX\{i\}} for \texttt{i=1,...,d}
using the barycentric Lagrange interpolation-based procedure. Lines 40-45 show
the computation of the columns of $N$ for each polynomial degree using the
built-in \texttt{null} function in MATLAB. The \texttt{null} function computes
the singular value decomposition of the input matrix in the background and sets
the output to the right singular vectors of the matrix that have singular values
close to machine epsilon.

\begin{figure}[h!]
\centering
\begin{minipage}{\columnwidth}
\begin{lstlisting}
function [N,Q,H,Qd,x,w,x1,C,DX,D1,v1]=ardivfreebfref(k,d)
%%%% function [N,Q,H,Qd,x,w,x1,C,DX,D1,v1]=ardivfreebfref(k,d)  
%%%% Inputs: k - Polynomial degree, d - Spatial dimension.
%%%% Outputs: N - Coeff. matrix, Qd - Div. free poly. basis, Q - Poly. basis, H - Upper Hess. matrix,
%%%% x - Int. pts, w - Int. wts, C - Div. matrix, DX - Der. matrix, D1 - 1d Der. matrix, v1 - Bary. wts.
  kp1=k+1; kp1d=(kp1)^d; kdim=nchoosek(k+d,d); kdimp=nchoosek(k-1+d,d);
  dkdim=d*kdim; kddim=dkdim-kdimp; wf=double(factorial(d)); hkdim=nchoosek(k+d-1,d-1);
%%%%--------------------Quadrature rule. 
  for i=1:d, [x1{i},w1{i}]=GaussJacobi(kp1,0,d-i); x1{i}=0.5*(1+x1{i}); w1{i}=(0.5)^(d-i+1)*w1{i}; end
  [xt{1:d}]=ndgrid(x1{:}); [wt{1:d}]=ndgrid(w1{:});
  xt=reshape(cat(d,xt{:}),[kp1d d]); wt=reshape(cat(d,wt{:}),[kp1d d]);
  x=[xt(:,1) xt(:,2:d).*cumprod(1-xt(:,1:d-1),2)]; w=prod(wt,2);clear xt wt;
%%%%--------------------Derivative matrices.
  for i=1:d, v1{i}=barywts(x1{i}); D1{i}=barydiff(x1{i},v1{i}); end
  DR{1}=D1{1};for i=2:d, DR{i}=speye(kp1); end
  for i=1:d,
    for j=2:d, if(j==i), DR{i}=kron(D1{i},DR{i}); else, DR{i}=kron(speye(kp1),DR{i}); end; end
  end
  DX{1}=DR{1}; sz=size(DX{1}); for i=2:d, DX{i}=zeros(sz); end; tt1=ones(kp1d,1);
  for j=2:d
    tt1=tt1-x(:,j-1);tt2=1./tt1;tt3=x(:,j)./tt1.^2;
    for i=1:j-1, DX{i}=DX{i}+DR{j}.*tt3; end; DX{j}=DX{j}+DR{j}.*tt2;
  end
%%%%--------------------Loop over polynomial degree.
  N=zeros(dkdim,kddim); N([1:kdim:dkdim],[1:d])=eye(d); C=zeros(kdimp,dkdim); DXQ=zeros(kp1d,d*hkdim);
  Q=zeros(kp1d,kdim); Q(:,1)=ones(kp1d,1); H=zeros(kdim,kdim-1); jdimpp=0; jdimp=1; jddimp=d; ct=1;
  if(nargout>3), Qd=zeros(d*kp1d,kddim); for i=1:d, Qd((i-1)*kp1d+1:i*kp1d,i)=1; end; end;  
  for j=1:k  
    jdim=nchoosek(j+d,d); djdim=d*jdim; jddim=djdim-jdimp; hj=jdim-jdimp; dhj=d*hj; hjddim=jddim-jddimp;
%%%%--------------------Arnoldi operations.    
    for i=1:d
      jjj=nchoosek(j-1+d-i,d-i);
      for jj=1:jjj
        q=x(:,i).*Q(:,jdimp-jjj+jj);  %% Multiply with x(:,i).
        [H(1:ct,ct),q]=mgs_with_reorth(Q,q,w,wf,ct,1);
        H(ct+1,ct)=sqrt(q'*(w.*q))*sqrt(wf); Q(:,ct+1)=q/H(ct+1,ct); ct=ct+1;   
      end
    end
%%%%--------------------Divergence matrix.
    for i=1:d, DXQ(:,(i-1)*hj+1:i*hj)=DX{i}*Q(:,jdimp+1:jdim); end
    ii=[]; for i=1:d, ii=[ii (i-1)*kdim+jdimp+1:(i-1)*kdim+jdim]; end    
    C(1:jdimp,ii)=mgs_with_reorth(Q(:,1:jdimp),DXQ(:,1:dhj),w,wf,jdimp,dhj);
%%%%--------------------Coefficient matrix.
    ii=[]; for i=1:d, ii=[ii (i-1)*kdim+1:(i-1)*kdim+jdim]; end
    N(ii,jddimp+1:jddim)=null([C(1:jdimp,ii);N(ii,1:jddimp)']);
%%%%--------------------Construct the basis.      
    if(nargout>3)
      for i=1:d, Qd((i-1)*kp1d+1:i*kp1d,jddimp+1:jddim)=...
                 Q(:,1:jdim)*N((i-1)*kdim+1:(i-1)*kdim+jdim,jddimp+1:jddim); end
    end    
    jdimpp=jdimp; jdimp=jdim; jddimp=jddim; 
  end
function [dd,q]=mgs_with_reorth(Q,f,w,wfact,n,nrhs)
  d=zeros(nrhs,n); q=f; 
  for i=1:n, t=Q(:,i)'*(w.*q)*wfact; q=q-Q(:,i).*t; d(:,i)=t.'; end
  for i=1:n, t=Q(:,i)'*(w.*q)*wfact; q=q-Q(:,i).*t; d(:,i)=d(:,i)+t.'; end; dd=d';
function [D]=barydiff(x,w,n)
  n=length(x); D=zeros(n,n); 
  for i=1:n, D(:,i)=(1/w(i))*(w(:)./(x(i)-x(:))); D(i,i)=-sum(D((1:n)~=i,i)); end; D=D';
function [w]=barywts(x)
  n=length(x); w=zeros(n,1); w(1)=1; 
  for j=2:n
    for k=1:j-1, w(k)=w(k)*(x(k)-x(j)); end; w(j)=prod(x(j)-x(1:j-1),'all');
  end; w=1./w;
\end{lstlisting}
\end{minipage}
\caption{\texttt{ardivfreebfref.m} - MATLAB code for \textbf{Step 1}.}
\label{fig:ardivfreebfref}
\end{figure}

\section {MATLAB implementation of step 2} \label{app:matlabstep2}

The MATLAB function \texttt{ardivfreebfgen} in figure \ref{fig:ardivfreebfgen}
shows our MATLAB implementation of \textbf{step 2}. The built-in function
\texttt{qr} computes the QR factorization of the input matrix and it returns the
orthonormal matrix as the first output argument. The column vectors of the
returned orthonormal matrix are orthonormal up to machine precision.

\begin{figure}[t]
\centering
\begin{minipage}{\columnwidth}
\begin{lstlisting}
function [Ne,Qde]=ardivfreebfgen(k,d,Xe,N,varargin)
%%%% function [Ne,Qde]=ardivfreebfgen(k,d,Xe,N,Q)  
%%%% Inputs: k - Polynomial degree, d - Spatial dimension. N - Ref. elem. coeff. matrix,
%%%% Q - Ref. elem. poly. basis, Xe - Gen. elem. node coord. matrix.
%%%% Outputs: Ne - Gen. elem. coeff. matrix, Qde - Gen. elem. div. free poly. basis.
  kp1=k+1; kp1d=(kp1)^d; kdim=nchoosek(k+d,d); kdimp=nchoosek(k-1+d,d); dkdim=d*kdim; kddim=dkdim-kdimp;
  Xet=Xe.'; Fe=Xet(:,2:d+1)-repmat(Xet(:,1),[1 d]); clear Xet;
  Feeye=kron(Fe,speye(kdim)); Ne=Feeye*N;  
%%%%--------------------Loop over polynomial degree.
  if(nargout>1), Qde=zeros(d*kp1d,kddim); end; jdimp=0; jddimp=0;Nbar=zeros(size(Ne));
  for j=0:k
    jdim=nchoosek(j+d,d);djdim=d*jdim;jddim=djdim-jdimp;hjddim=jddim-jddimp;
%%%%--------------------Linear combination of columns of N.
    ii=[]; for i=1:d, ii=[ii (i-1)*kdim+1:(i-1)*kdim+jdim]; end
    Nbar(1:djdim,1:jddim)=Ne(ii,1:jddim);
%%%%--------------------Orthonormalize the linear combinations.
    [Nbar(1:djdim,jddimp+1:jddim),~]=qr(Nbar(1:djdim,jddimp+1:jddim)-Nbar(1:djdim,1:jddimp)*(Nbar(1:djdim,1:jddimp)'*Nbar(1:djdim,jddimp+1:jddim)),0);
    T=Nbar(1:djdim,1:jddimp)'*Nbar(1:djdim,jddimp+1:jddim);tt=norm(T,'inf');
    if(tt > 1d-13) %% Needed for skewed simplices.
      fprintf('Reorthogonalizing... tt=%e\n',tt);
      [Nbar(1:djdim,jddimp+1:jddim),~]=qr(Nbar(1:djdim,jddimp+1:jddim)-Nbar(1:djdim,1:jddimp)*T,0);
    end
    Ne(ii,jddimp+1:jddim)=Nbar(1:djdim,jddimp+1:jddim);
%%%%--------------------Construct the basis.      
    if(nargout > 1)
      for i=1:d, Qde((i-1)*kp1d+1:i*kp1d,jddimp+1:jddim)=...
                 varargin{1}(:,1:jdim)*Ne((i-1)*kdim+1:(i-1)*kdim+jdim,jddimp+1:jddim); end
    end
    jdimp=jdim;jddimp=jddim;
  end    
\end{lstlisting}
\end{minipage}
\caption{\texttt{ardivfreebfgen.m} - MATLAB code for \textbf{step 2}.}
\label{fig:ardivfreebfgen}
\end{figure}

\section{MATLAB implementation of the evaluation algorithm} \label{app:matlabeval}

The MATLAB function \texttt{ardivfreebfeval} in figure
\ref{fig:ardivfreebfeval}, shows the evaluation algorithm. The mapped
coordinates of the points is specified using the $n_p\times d$ matrix
\texttt{s}, where $n_p$ is the number of points.


\begin{figure}[t]
\centering
\begin{minipage}{\columnwidth}
\begin{lstlisting}
function [Wd,W]=ardivfreebfeval(k,d,H,Ne,s)
%%%% Inputs: k - Poly. deg., d - Spatial dim., H - Upper Hess. matrix, Ne - Coeff. matrix, s - Eval. pts.
%%%% Outputs: Wd - Div. free basis at eval. pts, W - Poly. basis at eval. pts.
  kdim=nchoosek(k+d,d);dkdim=d*kdim;kdimp=nchoosek(k-1+d,d);kddim=dkdim-kdimp;M=size(s,1);
%%%%--------------------Loop over polynomial degree.      
  W=zeros(M,kdim);W(:,1)=ones(M,1);Wd=zeros(d*M,kddim);Wd(:,1:d)=kron(eye(d),W(:,1))*Ne([1:kdim:dkdim],[1:d]);
  ct=1;jdimp=1;jddimp=d;  
  for j=1:k
    jdim=nchoosek(j+d,d);jddim=d*jdim-jdimp;
    for i=1:d
      jjj=nchoosek(j-1+d-i,d-i);
      for jj=1:jjj
      w=s(:,i).*W(:,jdimp-jjj+jj);
      for ii=1:ct, w=w-H(ii,ct)*W(:,ii); end
      W(:,ct+1)=w/H(ct+1,ct);ct=ct+1;
      end
    end
    ii=[];for i=1:d, ii=[ii (i-1)*kdim+1:(i-1)*kdim+jdim]; end
    Wd(:,jddimp+1:jddim)=kron(speye(d),W(:,1:jdim))*Ne(ii,jddimp+1:jddim);
    jdimp=jdim;jddimp=jddim;
  end

\end{lstlisting}
\end{minipage}
\caption{\texttt{ardivfreebfeval.m} - MATLAB function to evaluate the
  divergence-free basis functions at points \texttt{s}.}
\label{fig:ardivfreebfeval}
\end{figure}

\end{document}